\documentclass[preprint,12pt]{elsarticle}

\usepackage{amssymb}
\usepackage{amsmath}
\usepackage{amsthm}
\usepackage{mathtools}
\usepackage{enumitem}
\usepackage{needspace}
\usepackage{tikz-cd}
\usepackage[colorlinks=true,linkcolor=blue,citecolor=blue,urlcolor=blue]{hyperref}
\usepackage[dvipsnames]{xcolor}

\journal{Journal of Algebra}

\theoremstyle{plain}
\newtheorem{theorem}{Theorem}[section]
\newtheorem{lemma}[theorem]{Lemma}
\newtheorem{proposition}[theorem]{Proposition}
\newtheorem{corollary}[theorem]{Corollary}
\theoremstyle{definition}
\newtheorem{definition}[theorem]{Definition}
\newtheorem{example}[theorem]{Example}
\newtheorem{remark}[theorem]{Remark}

\theoremstyle{plain}

\newcommand{\GV}{\operatorname{GV}}

\newcommand{\wMax}{w\operatorname{Max}}
\newcommand{\wSpec}{w\operatorname{Spec}}
\newcommand{\Ann}{\operatorname{Ann}}
\newcommand{\Hom}{\operatorname{Hom}}
\newcommand{\Tor}{\operatorname{Tor}}
\newcommand{\Ext}{\operatorname{Ext}}
\newcommand{\Spec}{\operatorname{Spec}}
\newcommand{\Max}{\operatorname{Max}}

\newcommand{\wfd}{w\text{-fd}}
\newcommand{\wwgldim}{w\text{-}w.\ gl.\ dim}
\newcommand{\wgl}{w.\ gl.\ dim}

\newcommand{\orcid}[1]{\,\href{https://orcid.org/#1}{\textsuperscript{[ORCID]}}\,}

\hypersetup{pdftitle={The w-operation and finiteness conditions in trivial ring extensions},
pdfauthor={Qing Li}}

\begin{document}

\begin{frontmatter}

\title{The \texorpdfstring{$w$}{w}-operation and finiteness conditions
in trivial ring extensions}

\author[swmu]{Qing Li\orcid{0000-0003-3490-7981}}
\ead{lqop80@163.com}
\address[swmu]{School of Mathematics, Southwest Minzu University,
Chengdu 610041, P. R. China}

\begin{abstract}
Let \(R\) be a commutative ring, \(M\) an \(R\)-module,
\(S=R\ltimes M\), and \(J=0\ltimes M\).
For a finitely generated ideal \(G\subseteq S\), set
\(I=\pi(G)\), where \(\pi:S\to R\) is the canonical projection.
Then \(G\) is a GV-ideal of \(S\) if and only if \(I\) is a
GV-ideal of \(R\) and the canonical multiplication map
\(M\to\Hom_R(I,M)\) is an isomorphism. Every GV-ideal \(I\) of \(R\) has \(IS\) GV in \(S\)
if and only if \(M\) is a \(w\)-module, equivalently
\(\wSpec(S)=\{\mathfrak p\ltimes M:\mathfrak p\in\wSpec(R)\}\).
For a \(w\)-module \(M\), an \(S\)-module \(X\) with \(JX=0\)
is of \(w\)-finitely presented type over \(S\)
if and only if it is so over \(R\) and \(M\otimes_R X\) is of
\(w\)-finite type. For an \(S\)-module \(X\), the conditions that the
\(R=S/J\)-module \(X/JX\) is of \(w\)-finitely presented type over
\(R\) and that \(\Tor_1^S(R,X)\) is of \(w\)-finite type are
sufficient for \(X\) to be of \(w\)-finitely presented type over
\(S\), and they are also necessary when \(R\) is \(w\)-coherent.
Furthermore, \(S\) is \(w\)-Noetherian if and only if \(R\) is
\(w\)-Noetherian and \(M\) is of \(w\)-finite type. If \(M\) is
of \(w\)-finite type, then \(S\) is \(w\)-coherent if and only if
\(R\) is \(w\)-coherent and \(M\) is of \(w\)-finitely presented
type. The \(w\)-weak global dimension of \(S\) is infinite when
\(M\) is a nonzero \(w\)-flat \(w\)-module of \(w\)-finite type.
For a domain \(D\) and a nonzero nonmaximal prime ideal \(P\), the
CPI-extension \(C(D,P)=D+PD_P\) is \(w\)-coherent if and only if
\(D/P\) is \(w\)-coherent and \(D_P\) is a valuation domain.
\end{abstract}

\begin{keyword}
$w$-operation \sep GV-ideal \sep trivial ring extension \sep pullback \sep
$w$-coherent ring \sep Milnor square \sep CPI-extension
\MSC[2020]{13A15; 13D05; 13D30; 13E05; 13F05}
\end{keyword}

\end{frontmatter}

\section{Introduction}
\label{sec:intro}

Throughout this paper, unless otherwise stated, all rings are nonzero
commutative rings with identity, and all modules are unital.
We use the $w$-module and finiteness terminology of Wang and Kim
\cite{wang-monograph1,wang-monograph2}; see also
\cite{wang-mccasland,w-fg-orig}. In particular, an $R$-module $M$
is called a $w$-module if it is GV-torsion-free and satisfies
$\Ext^1_R(R/I,M)=0$ for every GV-ideal $I$ of $R$.

For an $R$-module $M$, the trivial ring extension $R\ltimes M$
has underlying $R$-module $R\oplus M$ and multiplication
\[
(r,m)(r',m')=(rr',rm'+r'm).
\]
Throughout the paper, unless otherwise stated, we write
$S=R\ltimes M$ and $J=0\ltimes M$, and denote by
$\pi:S\to R$, $\pi(r,m)=r$, the canonical projection.
The quotient $S/J$ is $R$, and
$\mathfrak p\mapsto\mathfrak p\ltimes M$ identifies their ordinary
prime spectra. The $w$-operation also depends on the multiplication
maps $S\to\Hom_S(G,S)$ for finitely generated ideals $G\subseteq S$.
These maps involve the square-zero summand $M$, and the prime
correspondence need not identify the $w$-spectra. Thus idealizations
make it possible to study how the module structure affects the
$w$-operation while leaving the ordinary prime spectrum unchanged.
The first problem is to characterize the modules $M$ for which
GV-ideals lift and the $w$-maximal ideals correspond.
Classical idealization theory is developed in
\cite{anderson-winders,nagata}; see \cite{survey-te} for a survey.
Related finiteness questions for amalgamated algebras occur in
\cite{bak-kabbaj-mahdou,chhiti-jarrar-kabbaj-mahdou,danna,danna-finocchiaro-fontana,finocchiaro-prufer,kabbaj-louartiti-tamekkante},
and homological questions for trivial extensions in
\cite{adarbeh-kabbaj,mao-ding,gp-te-2025}.

Chang and Kim \cite{chang-kim} study GV-ideals, $w$-spectra and
the $w$-FF property, including $D\ltimes E$ for a domain $D$ and
a vector space $E$ over its quotient field.
Related correspondences for divisible $w$-modules over GCD domains
satisfying property (H) appear in
\cite[Propositions~4.4(3) and~4.5]{wang-zhou-chen-fractions};
the square-zero case of $D+\theta K[\theta]$ is treated in
\cite[Lemma~2.1]{zckz}.
Our first result applies to arbitrary $R$ and $M$.
For every finitely generated ideal $G\subseteq S$,
put $I=\pi(G)\subseteq R$.
Theorem~\ref{thm:gv-correspondence} gives
\[
 G\in\GV(S)\iff IS\in\GV(S)
 \iff I\in\GV(R)\text{ and }M\xrightarrow{\sim}\Hom_R(I,M).
\]
This includes nonhomogeneous ideals. Moreover, $M$ is a $w$-module
if and only if all GV-ideals lift, equivalently if the canonical
prime correspondence identifies the $w$-spectra
(Theorem~\ref{thm:wmax-correspondence}).
For arbitrary $M$, Theorem~\ref{thm:relative} describes the
relative family of GV-ideals and the closure of homogeneous ideals.
The examples illustrate the limitations of these criteria:
GV-torsion-freeness alone does
not ensure that $\eta_I^M$ is an isomorphism, and $I^2S\subseteq G$
cannot in general be strengthened to $IS\subseteq G$
(Examples~\ref{ex:torsionfree-not-wmodule} and~\ref{rem:sandwich-fails}).

The second problem is lifting finite presentations from $R=S/J$ to $S$.
The square-zero ideal can create additional relations, and finite
presentation at every $w$-maximal localization need not yield a module
of $w$-finitely presented type (Remark~\ref{rem:stalkwise-finiteness}).
The tensor and Tor terms below express this obstruction through global
$w$-finiteness conditions.
If $R$ is coherent, the following equivalence holds for every
$S$-module $X$ \cite[Theorem~4.4.2 and Corollary~4.4.3]{glaz}.
\[
 \begin{gathered}
 X\text{ is finitely presented over }S
 \\
 \Updownarrow
 \\
 \begin{cases}
 X/JX\text{ is finitely presented over }R,\\
 \Tor_1^S(R,X)\text{ is finitely generated over }R.
 \end{cases}
 \end{gathered}
\]
Assume now that $M$ is a $w$-module over $R$. We prove the following criteria.
\begin{enumerate}[label=(\roman*),leftmargin=*]
\item If $JX=0$, then $X$ is of $w$-finitely presented type over $S$
if and only if it is so over $R$ and $M\otimes_R X$ is of
$w$-finite type (Theorem~\ref{thm:inflated-fp}).
\item For an arbitrary $S$-module $X$, put $Q=X/JX$ and
$T=\ker(M\otimes_R Q\to JX)\cong\Tor_1^S(R,X)$.
If $Q$ is of $w$-finitely presented type over $R$ and $T$ is of
$w$-finite type, then $X$ is of $w$-finitely presented type over $S$.
The condition on $Q$ is necessary; that on $T$ is also necessary
when $R$ is $w$-coherent (Theorem~\ref{thm:w-tor-lifting}).
\end{enumerate}
Neither the tensor criterion nor the Tor criterion requires $M$ to
be of $w$-finite type.
Examples~\ref{ex:local-noetherian-not-global}--\ref{ex:tor-converse-fails}
show that finite generation at every $w$-maximal localization
does not imply $w$-finite type, that the tensor condition cannot be omitted, and
that the Tor converse can fail without $w$-coherence of the base ring.

Section~\ref{sec:transfer} also proves that $S$ is $w$-Noetherian
if and only if $R$ is $w$-Noetherian and $M$ is of $w$-finite type.
$w$-coherence descends from $S$ to $R$ without a finite-type assumption
on $M$. If $M$ is of $w$-finite type, or contains a finite family
with zero common annihilator, $S$ is $w$-coherent if and only if
$R$ is $w$-coherent and $M$ is of $w$-finitely presented type.
Examples with divisible modules \cite{adarbeh-kabbaj,wang-kim-relative-if}
show that $w$-coherence of $S$ alone does not imply that $M$ is of
$w$-finite type. The $w$-local Gorenstein statements apply
the classical Reiten--Fossum theorem after the global
$w$-Noetherian condition has been checked.

Sections~\ref{sec:patching} and~\ref{sec:cpi} give square-zero and
pullback applications. Classical pullback theory includes $D+M$
rings and CPI-extensions
\cite{bastida-gilmer,boisen-sheldon,boisen-sheldon-preprufer,brewer-rutter,dobbs-papick},
and homological and coherence transfer
\cite{costa-kabbaj,fontana-gabelli,gabelli-houston}.
For related $w$-theoretic and Pr\"ufer-type results, see
\cite{chang-kim-prufer,lin-wang-wphi,wang-wR,wang-omega,xie-wang-tian,xie-tang}.
We distinguish the square-zero kernel from the conductor and apply
the classical Tor formula for tensor powers \cite{fgr} to weak global
dimension. The CPI results specialize established Milnor-square
criteria to $C(D,P)=D+PD_P$, including $P=0$ and maximal $P$.

\section{Preliminaries}
\label{sec:prelim}

We use Wang and Kim's definitions
\cite[\S\,6.1 and \S\,6.3]{wang-monograph2}.
A finitely generated ideal $I$ of $R$ is a \emph{GV-ideal} if the
canonical multiplication map
\[
R\longrightarrow\Hom_R(I,R),\qquad r\longmapsto(a\mapsto ra)
\]
is an isomorphism. We denote these ideals by $\GV(R)$.
They are faithful and are closed under products.
For an $R$-module $X$, put
\[
\operatorname{tor}_{\GV}(X)
=\{x\in X:Ix=0\text{ for some }I\in\GV(R)\}.
\]
The module $X$ is \emph{GV-torsion} if
$\operatorname{tor}_{\GV}(X)=X$, and \emph{GV-torsion-free} if
$\operatorname{tor}_{\GV}(X)=0$.
Faithfulness of GV-ideals implies $\operatorname{tor}_{\GV}(R)=0$.

An $R$-module $X$ is called a \emph{$w$-module} if it is GV-torsion-free and satisfies
\[
\Ext^1_R(R/I,X)=0\qquad(I\in\GV(R))
\quad\text{\cite[Definition~6.1.12]{wang-monograph2}}.
\]
By applying
$\Hom_R(-,X)$ to $0\to I\to R\to R/I\to0$, it is equivalent
to requiring that $X\to\Hom_R(I,X)$ be an isomorphism for every
$I\in\GV(R)$. Thus every $R$-linear map $I\to X$ extends uniquely
to $R$.
Flat modules are $w$-modules \cite[Theorem~6.7.35]{wang-monograph2}.

For GV-torsion-free $X$, its \emph{$w$-envelope} is
\[
X_w=\{x\in E(X):Ix\subseteq X\text{ for some }I\in\GV(R)\},
\]
where $E(X)$ is its injective hull. Then $X$ is a $w$-module if and
only if $X=X_w$ \cite[Theorem~6.1.15]{wang-monograph2}.
For an ideal $A$ of $R$, we use the closure inside $R$,
\[
A_w=\{r\in R:Ir\subseteq A\text{ for some }I\in\GV(R)\}.
\]
In particular, $R_w=R$.

A prime is a \emph{$w$-prime} if it is $w$-closed. We write
$\wSpec(R)$ for the set of such primes and $\wMax(R)$ for the
maximal proper $w$-ideals, which are prime. Every proper $w$-ideal
is contained in a $w$-maximal ideal. A prime is $w$-closed exactly
when it contains no GV-ideal. A finitely generated ideal $I$ is a
GV-ideal exactly when $I_w=R$.
The localization criterion is \cite[Theorem~6.2.13]{wang-monograph2}
\[
X\text{ is GV-torsion}\iff
X_{\mathfrak m}=0\quad\text{for all }\mathfrak m\in\wMax(R).
\]
A sequence of $R$-modules is \emph{$w$-exact} if its localization
at every $\mathfrak m\in\wMax(R)$ is exact
\cite[Definition~6.3.1]{wang-monograph2}. In particular, a sequence
\[
0\longrightarrow X'\xrightarrow{f}X\xrightarrow{g}X''
\longrightarrow0
\]
is $w$-exact if and only if
\[
0\longrightarrow X'_{\mathfrak m}\xrightarrow{f_{\mathfrak m}}
X_{\mathfrak m}\xrightarrow{g_{\mathfrak m}}X''_{\mathfrak m}
\longrightarrow0
\]
is exact for every $\mathfrak m\in\wMax(R)$.
For a submodule $Y$ of a $w$-module $X$, its closure in $X$
satisfies
\[
Y_w=\{x\in X:x/1\in Y_{\mathfrak m}
\ \forall\,\mathfrak m\in\wMax(R)\},
\]
by \cite[Theorems~6.2.14 and~6.2.15]{wang-monograph2}.
A module $X$ is \emph{of $w$-finite type} if there is a finitely
generated submodule $B\subseteq X$ with $X/B$ GV-torsion.
Equivalently, there exist a finitely generated free $R$-module $F$
and an $R$-linear map $f:F\to X$ such that
\[
f_{\mathfrak m}:F_{\mathfrak m}\twoheadrightarrow X_{\mathfrak m}
\qquad\bigl(\mathfrak m\in\wMax(R)\bigr);
\]
that is, $f$ is a $w$-epimorphism.
A GV-torsion-free module $X$ is of $w$-finite type if and only if
\[
B_w=X_w\qquad\text{for some finitely generated }B\subseteq X.
\]
A module $X$ is \emph{of $w$-finitely presented type} if there is a $w$-exact
sequence $F_1\to F_0\to X\to0$ with $F_0,F_1$ finitely generated
free. These are the notions called ``finite type'' and ``finitely
presented type'' in \cite{wang-monograph2}.
A module $X$ is \emph{finitely presented} if there is an exact sequence
\[
F_1\longrightarrow F_0\longrightarrow X\longrightarrow0
\]
with $F_0,F_1$ finitely generated free.

An $R$-module $X$ of $w$-finite type is called \emph{$w$-coherent}
if every submodule of $X$ of $w$-finite type is of $w$-finitely
presented type \cite[Definition~6.11.1]{wang-monograph2}.

A ring is a \emph{DW-ring} if its $w$-operation is the identity,
equivalently if $\GV(R)=\{R\}$ or $\wMax(R)=\Max(R)$
\cite[Theorem~6.2.20]{wang-monograph2}.
On such a ring the two finite-type notions above coincide with
finite generation and finite presentation.
Zero-dimensional rings, one-dimensional domains and Pr\"ufer
domains are DW-rings. A local ring need not be DW. In
$k[[x,y]]$, the maximal ideal $(x,y)$ is a GV-ideal.

We use the following criteria.

\Needspace{6\baselineskip}
\begin{proposition}
\label{prop:finite-type-tools}
Let $A$ be a ring.
\begin{enumerate}[label=(\arabic*),leftmargin=*]
\item If $0\to K\to F\to X\to0$ is $w$-exact and $F$ is
finitely generated free, then $X$ is of $w$-finitely presented type
if and only if $K$ is of $w$-finite type.
\item If $0\to K\to B\to X\to0$ is $w$-exact, $B$ is of
$w$-finite type and $X$ is of $w$-finitely presented type, then
$K$ is of $w$-finite type.
\item If $B$ is $w$-coherent and $K\subseteq B$ is of $w$-finite type,
then $K$ and $B/K$ are $w$-coherent.
\item A ring $A$ is $w$-coherent if and only if every finitely
generated ideal is of $w$-finitely presented type. Over such a ring,
every module of $w$-finitely presented type is $w$-coherent.
\end{enumerate}
\end{proposition}

\begin{proof}
See \cite[Theorems~6.3.22, 6.3.25 and~6.11.5,
and Proposition~6.11.3]{wang-monograph2}.
\end{proof}

\begin{remark}
\label{rem:stalkwise-finiteness}
For an $A$-module $N$, the condition that $N_{\mathfrak m}$ is
finitely presented over $A_{\mathfrak m}$ for every
$\mathfrak m\in\wMax(A)$ does not imply that $N$ is of
$w$-finitely presented type, even if $N$ is cyclic. Let $A=\prod_{i\ge1}k$, let
$I=\bigoplus_{i\ge1}k$, and put $N=A/I$. Every prime localization
of $A$ is a field. If a prime omits some coordinate idempotent
$e_i$, then $I$ localizes to the whole ring; otherwise all $e_i$
vanish in the localization and $I$ localizes to zero. Hence every
localization of the cyclic module $N$ is zero or a field, and is
finitely presented. But $I$ is not finitely generated, so $N=A/I$
is not finitely presented. Since $A$ is zero-dimensional and hence
a DW-ring, $N$ is not of $w$-finitely presented type.
\end{remark}

The ring $R$ is \emph{$w$-Noetherian} if every ideal is of
$w$-finite type, equivalently if every $w$-ideal is of $w$-finite
type. On domains this is the strong Mori condition
\cite{wang-mccasland}. In particular,
\begin{equation}\label{eq:w-noetherian}
R\text{ is }w\text{-Noetherian}\implies
R_{\mathfrak m}\text{ is Noetherian for all }\mathfrak m\in\wMax(R).
\end{equation}
Every module of $w$-finite type over
a $w$-Noetherian ring is of $w$-finitely presented type and is
$w$-Noetherian; finite direct sums and submodules of
$w$-Noetherian modules are $w$-Noetherian
\cite[\S\,6.8]{wang-monograph2}. Further homological descriptions
appear in \cite{w-hom-mod}.

A module is $w$-flat if its localization at every $w$-maximal
ideal is flat. We use the dimension formulas
\[
\begin{aligned}
\wfd_R X&=\sup_{\mathfrak m\in\wMax(R)}
\operatorname{fd}_{R_{\mathfrak m}}X_{\mathfrak m},\\
\wwgldim(R)&=\sup_{\mathfrak m\in\wMax(R)}\wgl(R_{\mathfrak m})
\end{aligned}
\]
from \cite{wang-monograph2,w-wgd}.

Let $(A,\mathfrak m)$ be a $d$-dimensional Cohen--Macaulay
Noetherian local ring, and put $k=A/\mathfrak m$.
A \emph{canonical module} of $A$ is a finitely generated
$A$-module $\omega_A$ satisfying \cite[\S\,3.3]{bruns-herzog}
\[
\operatorname{depth}_A\omega_A=d,\qquad
\operatorname{id}_A\omega_A<\infty,\qquad
\dim_k\Ext_A^d(k,\omega_A)=1.
\]

\begin{definition}
\label{def:w-cm-gorenstein}
We use the following $w$-local terminology.
\begin{enumerate}[label=(\roman*),leftmargin=*]
\item A ring is $w$-Cohen--Macaulay if it is $w$-Noetherian and
all its localizations at $w$-maximal ideals are Cohen--Macaulay.
\item A ring is $w$-Gorenstein if it is $w$-Noetherian and all
these localizations are Gorenstein.
\item A $w$-canonical module over a $w$-Cohen--Macaulay ring $R$
is a $w$-module $C$ such that $C_{\mathfrak m}$ is a canonical
module of $R_{\mathfrak m}$ for every $\mathfrak m\in\wMax(R)$.
\end{enumerate}
\end{definition}

\begin{remark}\label{rem:wGor-terminology}
These definitions impose no uniform bound on the local injective
dimensions. On a DW-ring these definitions reduce to the usual
Noetherian Cohen--Macaulay and Gorenstein notions.
\end{remark}

We use a cartesian square in the form \cite{milnor}
\begin{equation}\label{eq:milnor-square}
\begin{tikzcd}
R \ar[r,hook] \ar[d,two heads] & T \ar[d,two heads]\\
D \ar[r,hook] & F.
\end{tikzcd}
\end{equation}
Here $R=T\times_F D$ and $\mathfrak c=\ker(T\to F)$ is a
common ideal of $R$ and $T$, with $R/\mathfrak c\cong D$.
We call it the common kernel. It is contained in the full conductor
$(R:T)=\{r\in R:rT\subseteq R\}$.

The domain Milnor squares of \cite[Chapter~8]{wang-monograph2}
have a nonzero maximal ideal $\mathfrak c$ of a domain $T$, a
field $F=T/\mathfrak c$, and a proper subring $D\subsetneq F$.
Their conductor equals $\mathfrak c$. They are of type I when
$F=\operatorname{Frac}(D)$ and of type II otherwise.

Let $S=R\ltimes M$ and $J=0\ltimes M$. We identify $R$ with
$R\ltimes0$ and use the following standard facts
\cite[Theorems~3.2, 3.3, 3.5 and~4.1]{anderson-winders}.
\begin{enumerate}[label=(\arabic*),leftmargin=*]
\item Every prime contains $J$. The map
$\mathfrak p\mapsto\mathfrak p\ltimes M$ is a homeomorphism
$\Spec(R)\to\Spec(S)$ and takes maximal ideals to maximal ideals.
\item A homogeneous subset $I\oplus N$ is an ideal if and only
if $I$ is an ideal of $R$, $N$ is a submodule of $M$, and
$IM\subseteq N$. General ideals need not be homogeneous.
\item For every prime $\mathfrak p$ of $R$,
\[
S_{\mathfrak p\ltimes M}\cong R_{\mathfrak p}\ltimes M_{\mathfrak p}
\cong S\otimes_R R_{\mathfrak p}.
\]
This identifies the corresponding localizations of every $S$-module.
\end{enumerate}

For a finitely generated module $N$, the classical criteria are
\[
\begin{aligned}
A\ltimes N\text{ is Noetherian}
&\iff A\text{ is Noetherian},\\
A\ltimes N\text{ is coherent}
&\iff A\text{ is coherent and }N\text{ is finitely presented}.
\end{aligned}
\]
The first is \cite[Theorem~4.8]{anderson-winders}; the second follows
from \cite[p.~55 and Theorem~4.1.5, together with the observation
following Theorem~4.4.4]{glaz}.

\begin{lemma}
\label{lem:classical-reiten-fossum}
Let $A$ be a Noetherian local ring and let $N$ be a nonzero
finitely generated $A$-module. Then $A\ltimes N$ is Gorenstein
if and only if $A$ is Cohen--Macaulay and $N$ is a canonical
$A$-module.
\end{lemma}

\begin{proof}
This is the Reiten--Fossum theorem \cite{fossum,reiten};
see also \cite[Theorem~5.6]{fgr}.
\end{proof}

When $N=0$, the ring $A\ltimes N\cong A$ is Gorenstein if and only if
$A$ is Gorenstein. We give the homological calculation
needed for weak global dimension in Section~\ref{sec:patching}.

\section{The \texorpdfstring{$w$}{w}-operation on \texorpdfstring{$R\ltimes M$}{R ltimes M}}
\label{sec:internal}

Let \(M\) be an arbitrary \(R\)-module. Put \(S=R\ltimes M\), let
\(\pi:S\to R\) be projection, and write \(J=0\ltimes M\).
For an \(R\)-module \(X\) and a finitely generated ideal
\(I\subseteq R\), put
\[
 \eta_I^X:X\longrightarrow\Hom_R(I,X),\qquad x\longmapsto(a\mapsto ax).
\]
Extension along \(R\to S\) gives \(IS=I\ltimes IM\).

\begin{lemma}[{\cite[Proposition~6.1.9(3), (4)]{wang-monograph2}}]\label{lem:permanence}
Let \(A\) be a ring. Products of GV-ideals of \(A\) are GV-ideals of \(A\).
If \(U\in\GV(A)\) and \(V\) is a finitely generated ideal of \(A\)
containing \(U\), then \(V\in\GV(A)\).
\end{lemma}

\begin{lemma}\label{lem:powers}
Let \(G\) be an ideal of \(S\), and set \(I=\pi(G)\).
Then
\[
 G^2\subseteq IS,\qquad (IS)^2=I^2S\subseteq G.
\]
If \(G\) is finitely generated, then
\[
 G\in\GV(S)\quad\Longleftrightarrow\quad IS\in\GV(S).
\]
\end{lemma}
\begin{proof}
For \((a,x),(b,y)\in G\), we have \(a,b\in I\), and hence
\[
 (a,x)(b,y)=(ab,ay+bx)\in I^2\ltimes IM\subseteq IS.
\]
Thus \(G^2\subseteq IS\). Conversely, for any \(a,b\in I\),
choose lifts \((a,x),(b,y)\in G\). Since \(G\) is an ideal,
\[
 (a,x)(b,0)-(b,y)(0,x)=(ab,0)\in G.
\]
It follows that \(I^2\ltimes 0\subseteq G\). Multiplying by
elements of \(0\ltimes M\) gives \(0\ltimes I^2M\subseteq G\),
so \(I^2S=I^2\ltimes I^2M\subseteq G\). Also,
\((IS)^2=(I\ltimes IM)^2=I^2\ltimes I^2M=I^2S\).

Now suppose that \(G\) is finitely generated. Then \(I\) and
\(IS\) are finitely generated. If \(G\in\GV(S)\), then
\(G^2\in\GV(S)\); since \(G^2\subseteq IS\),
Lemma~\ref{lem:permanence} gives \(IS\in\GV(S)\).
Conversely, if \(IS\in\GV(S)\), then \((IS)^2\in\GV(S)\).
The inclusion \((IS)^2\subseteq G\) and finite generation of \(G\)
give \(G\in\GV(S)\), again by Lemma~\ref{lem:permanence}.
\end{proof}

\begin{theorem}\label{thm:projection}
Let \(G\in\GV(S)\) and put \(I=\pi(G)\). Then
\[
 I\in\GV(R),\qquad
 \eta_I^M:M\longrightarrow\Hom_R(I,M)
 \text{ is an isomorphism}.
\]
\end{theorem}
\begin{proof}
Generators of \(G\) project to generators of \(I\).
We first prove
\(\Ann_R(I)=0\) and \((0:_M I)=0\).
If \(Im=0\), then \((0,m)G=0\), so \(m=0\).
If \(rI=0\), then \((0,rm)G=0\) for every \(m\in M\);
thus \(rM=0\). Hence \((r,0)G=0\), and \(r=0\).
We used faithfulness of the GV-ideal \(G\).

For \(h\in\Hom_R(I,M)\), the map
\[
 F_h:G\to S,\qquad (a,x)\mapsto(0,h(a))
\]
is \(S\)-linear: for \((b,y)\in S\),
\(F_h((b,y)(a,x))=(0,bh(a))=(b,y)F_h(a,x)\).
It is multiplication by some \((r,n)\in S\).
First components give \(rI=0\), so \(r=0\);
second components give \(h(a)=an\).
Uniqueness of \(n\) follows from \((0:_M I)=0\).
Therefore \(\eta_I^M\) is an isomorphism.

Now let \(g\in\Hom_R(I,R)\). For each \(m\in M\), the map
\(a\mapsto g(a)m\) has a unique expression \(a\mapsto a\,u_g(m)\).
The inverse of \(\eta_I^M\) is \(R\)-linear, so \(u_g:M\to M\)
is \(R\)-linear and satisfies
\begin{equation}\label{eq:compat}
 a\,u_g(m)=g(a)m\qquad(a\in I,\ m\in M).
\end{equation}
Define \(F_g:G\to S\) by \(F_g(a,x)=(g(a),u_g(x))\).
For \((b,y)\in S\), equation~\eqref{eq:compat} gives
\[
\begin{aligned}
 F_g((b,y)(a,x))
   &=(bg(a),\,b\,u_g(x)+a\,u_g(y))\\
   &=(bg(a),\,b\,u_g(x)+g(a)y)\\
   &=(b,y)F_g(a,x).
\end{aligned}
\]
Thus \(F_g\) is \(S\)-linear and is multiplication by some
\((c,z)\in S\). First components give \(g(a)=ca\).
The scalar is unique since \(\Ann_R(I)=0\).
Hence \(\eta_I^R\) is an isomorphism and \(I\in\GV(R)\).
\end{proof}

\begin{theorem}\label{thm:gv-correspondence}
Let \(G\subseteq S\) be a finitely generated ideal and \(I=\pi(G)\).
The following are equivalent:
\begin{enumerate}[label=(\roman*)]
\item \(G\in\GV(S)\).
\item \(IS=I\ltimes IM\in\GV(S)\).
\item \(I\in\GV(R)\) and \(\eta_I^M\) is an isomorphism.
\end{enumerate}
\end{theorem}
\begin{proof}
Lemma~\ref{lem:powers} gives (i)\(\Leftrightarrow\)(ii), and
Theorem~\ref{thm:projection} gives (i)\(\Rightarrow\)(iii).
Under (iii), take \(f\in\Hom_S(IS,S)\) and write
\(f(a,0)=(\alpha(a),\beta(a))\) for \(a\in I\).
The maps \(\alpha:I\to R\) and \(\beta:I\to M\) are \(R\)-linear,
so \(\alpha(a)=ra\), \(\beta(a)=an\) for unique \(r\in R,n\in M\).
The elements \((a,0)\) generate \(IS\) over \(S\).
Thus \(f\) is multiplication by \((r,n)\) on all of \(IS\),
proving (ii).
\end{proof}

The implication (i)$\Rightarrow$(ii) also follows from
\cite[Theorem~6.9.2(7)$\Rightarrow$(6)]{wang-monograph2}, applied to
the canonical inclusion $R\hookrightarrow S$, since
$S=R\oplus M$ as an $R$-module.

\begin{theorem}\label{thm:wmax-correspondence}
For arbitrary \(M\), the following conditions are equivalent:
\begin{enumerate}[label=(\roman*)]
\item \(M\) is a \(w\)-module over \(R\).
\item \(IS\in\GV(S)\) for every \(I\in\GV(R)\).
\item For every finitely generated ideal \(G\subseteq S\),
\(G\in\GV(S)\) if and only if \(\pi(G)\in\GV(R)\).
\item \(\wSpec(S)=\{\mathfrak p\ltimes M:\mathfrak p\in\wSpec(R)\}\).
\item \(\wMax(S)=\{\mathfrak m\ltimes M:\mathfrak m\in\wMax(R)\}\).
\end{enumerate}
\end{theorem}
\begin{proof}
The equivalence of (i)--(iii) follows from
Theorem~\ref{thm:gv-correspondence}, since \(M\) is a \(w\)-module
if and only if \(\eta_I^M:M\to\Hom_R(I,M)\) is an isomorphism
for every \(I\in\GV(R)\).
For every \(M\), if \(\mathfrak p\in\wSpec(R)\), then
\(\mathfrak p\ltimes M\in\wSpec(S)\).
Indeed, if \(G\in\GV(S)\) and \(G\subseteq\mathfrak p\ltimes M\),
then Theorem~\ref{thm:projection} gives
\[
 \pi(G)\in\GV(R),\qquad \pi(G)\subseteq\mathfrak p,
\]
contradicting \(\mathfrak p\in\wSpec(R)\).
Under (ii), if \(\mathfrak p\ltimes M\in\wSpec(S)\), then
\(\mathfrak p\in\wSpec(R)\). Otherwise, there exists
\(I\in\GV(R)\) with \(I\subseteq\mathfrak p\), and (ii) gives
\(IS\in\GV(S)\) with \(IS\subseteq\mathfrak p\ltimes M\),
contradicting \(\mathfrak p\ltimes M\in\wSpec(S)\).
This proves (iv).

Assume (iv). If \(I\in\GV(R)\) but \(IS\notin\GV(S)\),
the proper \(w\)-ideal \((IS)_w\) is contained in some
\(\mathfrak q\in\wMax(S)\).
By (iv), \(\mathfrak q=\mathfrak p\ltimes M\) for a
\(w\)-prime \(\mathfrak p\) of \(R\).
Then \(I\subseteq\mathfrak p\), a contradiction.
This proves (ii). The same argument proves (v)\(\Rightarrow\)(ii).
Finally, (iv) implies (v), since the prime correspondence
preserves inclusion and \(w\)-maximal ideals are the maximal
elements among \(w\)-primes.
\end{proof}

For an arbitrary \(R\)-module \(M\), put
\[
 \mathcal G_M=\{I\in\GV(R):\eta_I^M\text{ is an isomorphism}\}
             =\{\pi(G):G\in\GV(S)\}.
\]
This family contains \(R\), is closed under products, and satisfies
\[
 I\in\mathcal G_M,\quad I\subseteq L\subseteq R,\quad
 L\text{ a finitely generated ideal}
 \quad\Longrightarrow\quad L\in\mathcal G_M.
\]
These assertions follow by extension to \(S\), using
Lemma~\ref{lem:permanence} and Theorem~\ref{thm:gv-correspondence}.
For a submodule \(Y\subseteq X\), define
\[
 Y_{c_M,X}
 =\{x\in X:Ix\subseteq Y\text{ for some }I\in\mathcal G_M\}.
\]

\Needspace{6\baselineskip}
\begin{theorem}\label{thm:relative}
Let \(M\) be an arbitrary \(R\)-module.
\begin{enumerate}[label=(\arabic*),leftmargin=*,nosep]
\item The \(w\)-spectrum of \(S\) is
\[
 \wSpec(S)=\{\mathfrak p\ltimes M:\mathfrak p\in\Spec(R),\
 I\nsubseteq\mathfrak p\text{ for all }I\in\mathcal G_M\}.
\]
\item If \(A\ltimes N\) is a homogeneous ideal, with \(AM\subseteq N\),
then its closure inside \(S\) is
\[
 (A\ltimes N)_{w_S}=A_{c_M,R}\ltimes N_{c_M,M}.
\]
In particular \(0\ltimes M\) is a \(w\)-ideal for every \(M\).
\end{enumerate}
\end{theorem}
\begin{proof}
Part (1) follows from Theorem~\ref{thm:gv-correspondence}.
Every GV-ideal of $S$ contains one of the form $IS$ with
$I\in\mathcal G_M$. Indeed, if $G\in\GV(S)$ projects to $I$,
then $I^2\in\mathcal G_M$ and $I^2S\subseteq G$.
Hence \((r,m)\in(A\ltimes N)_{w_S}\) precisely when some such
\(I\) satisfies
\[
 Ir\subseteq A,\qquad Im\subseteq N,\qquad rIM\subseteq N.
\]
The third condition follows from the first and \(AM\subseteq N\).
If $I_1r\subseteq A$ and $I_2m\subseteq N$ for
$I_1,I_2\in\mathcal G_M$, then $I_1I_2\in\mathcal G_M$
satisfies both conditions.
Finally, every \(I\in\mathcal G_M\) is faithful on \(R\),
so \(0_{c_M,R}=0\).
\end{proof}

If \(M\) is a \(w\)-module, then \(\mathcal G_M=\GV(R)\).
Thus the last formula recovers the componentwise \(w\)-closure
formula. For other modules, \(\mathcal G_M\) can be a proper subset of \(\GV(R)\).

\begin{corollary}\label{lem:lift-gv}
If \(M\) is a \(w\)-module, every \(I\in\GV(R)\) lifts to
the GV-ideal \(IS\) of \(S\).
\end{corollary}
\begin{proof}
Let \(I\in\GV(R)\). Since \(M\) is a \(w\)-module, the canonical map
\[
 \eta_I^M:M\longrightarrow\Hom_R(I,M)
\]
is an isomorphism.
Applying Theorem~\ref{thm:gv-correspondence} to the finitely generated
ideal \(IS\), we conclude that \(IS\in\GV(S)\).
\end{proof}

\begin{corollary}\label{cor:gv-projection}
Suppose \(M\) is a \(w\)-module.
\begin{enumerate}[label=(\arabic*),leftmargin=*,nosep]
\item For any prime \(\mathfrak p\)
of \(R\), a GV-ideal of \(S\) is contained in
\(\mathfrak p\ltimes M\) if and only if \(\mathfrak p\)
contains a GV-ideal of \(R\).
\item If \(G\in\GV(S)\) projects to \(I\), then \(G\) is
\(w\)-dense in \(I\ltimes M\).
\end{enumerate}
\end{corollary}
\begin{proof}
The first assertion follows by projection and lifting.
For the second, \(G_w=S\), since \(G\) is a GV-ideal;
in particular every element of \(I\ltimes M\) is in its closure.
Equivalently, \((I\ltimes M)/G\) is GV-torsion.
\end{proof}

\Needspace{6\baselineskip}
\begin{corollary}\label{cor:w-homeo}
Let \(M\) be a \(w\)-module.
\begin{enumerate}[label=(\arabic*),leftmargin=*,nosep]
\item The map
\(\mathfrak p\mapsto\mathfrak p\ltimes M\) gives an
order-preserving homeomorphism \(\wSpec(R)\to\wSpec(S)\), and
\begin{equation}\label{eq:w-dim-equal}
 w\text{-}\dim(S)=w\text{-}\dim(R).
\end{equation}
\item At every \(\mathfrak m\in\wMax(R)\),
\begin{equation}\label{eq:localization-formula}
 S_{\mathfrak m\ltimes M}\cong R_{\mathfrak m}\ltimes M_{\mathfrak m}.
\end{equation}
\end{enumerate}
\end{corollary}
\begin{proof}
Theorem~\ref{thm:wmax-correspondence} identifies the two
\(w\)-spectra inside the usual prime \mbox{spectrum} homeomorphism.
For any finite chain
\(\mathfrak p_0\subsetneq\cdots\subsetneq\mathfrak p_n\)
in \(\wSpec(R)\), the map in (1) gives a chain
\(\mathfrak p_0\ltimes M\subsetneq\cdots\subsetneq
\mathfrak p_n\ltimes M\) in \(\wSpec(S)\).
Conversely, contracting each prime in a chain in \(\wSpec(S)\)
recovers a chain in \(\wSpec(R)\). Thus the two spectra have the
same supremum of chain lengths, and
\(w\text{-}\dim(S)=w\text{-}\dim(R)\).
The localization formula is the usual one for idealization.
\end{proof}

\begin{corollary}\label{cor:w-ideals-lift}
If \(M\) is a \(w\)-module, then for every ideal \(A\subseteq R\),
\((A\ltimes M)_w=A_w\ltimes M\).
In particular, \(A\ltimes M\) is a \(w\)-ideal if and only if
\(A\) is a \(w\)-ideal.
\end{corollary}
\begin{proof}
Apply Theorem~\ref{thm:relative}; for a \(w\)-module,
\(\mathcal G_M=\GV(R)\).
\end{proof}

\begin{corollary}\label{cor:w-closure}
If \(M\) is a \(w\)-module and \(A\ltimes N\) is a homogeneous
ideal of \(S\), then
\begin{equation}\label{eq:w-closure-formula}
 (A\ltimes N)_w=A_w\ltimes N(M)_w,
\end{equation}
where \(N(M)_w=\{m\in M:Im\subseteq N\text{ for some }I\in\GV(R)\}\).
Thus \(A\ltimes N\) is a \(w\)-ideal if and only if
\(A=A_w\) and \(N=N(M)_w\).
\end{corollary}
\begin{proof}
This is Theorem~\ref{thm:relative} with
\(\mathcal G_M=\GV(R)\).
\end{proof}

\begin{remark}\label{rem:chang-kim-spectra}
For a domain \(D\) with quotient field \(K\), every
\(K\)-vector space \(E\) is a \(w\)-module over \(D\).
The preceding results recover Chang and Kim's GV-ideal, closure
and spectral calculations in that case
\cite[Proposition~2.2(4)--(6) and Corollary~2.3]{chang-kim}.
For $E=K$, these are also the square-zero case of
\cite[Lemma~2.1]{zckz}. Another established special case concerns a
GCD domain $D$ and a divisible $w$-module $H$ such that every
nonzero nonunit of $D$ annihilates a nonzero element of $H$.
Under these assumptions, the GV-ideal correspondence, the formula for $(A\ltimes H)_w$,
and the maximal $w$-ideal correspondence occur in
\cite[Propositions~4.4(3) and~4.5]{wang-zhou-chen-fractions}.
The $w$-finite type criterion for ideals of $D\ltimes E$ with
nonzero projection is given in \cite[Proposition~2.2(6)]{chang-kim}.
\end{remark}

\begin{example}\label{rem:sandwich-fails}
Let \(R=k[[x,y]]\), \(M=R\), and
\(G=S(x,1)+S(y,1)\). Its projection \(I=\pi(G)=(x,y)\) is a GV-ideal;
\(M=R\) is a \(w\)-module, so \(G\in\GV(S)\).
Nevertheless \((x,0)\notin G\). A representation would give
\(xa+yb=x\) and \(a+b\in(x,y)\), whereas the first equality
modulo \((x,y)^2\) forces the constant terms of \(a,b\) to be
\(1,0\). Thus \(IS\subseteq G\) cannot replace
\(I^2S\subseteq G\).
\end{example}

\begin{example}\label{ex:torsionfree-not-wmodule}
Let \(R=k[x,y]\) and \(M=I=(x,y)\).
Then \(M\) is torsion-free, hence GV-torsion-free.
The ideal \(I\) is a GV-ideal. However, for $m\in M$,
\[
 \eta_I^M(m)=\operatorname{id}_I
 \;\Longrightarrow\; xm=x
 \;\Longrightarrow\; m=1\notin M.
\]
Thus \(\eta_I^M\) is not
surjective and \(IS=I\ltimes I^2\) is not a GV-ideal.
Explicitly, the \(S\)-linear map \((a,u)\mapsto(0,a)\) on
\(IS\) is not a multiplication map. Thus $M$ is GV-torsion-free but
$\Ext^1_R(R/I,M)\cong\operatorname{Coker}\eta_I^M\ne0$, so $M$ is
not a $w$-module.
\end{example}

\section{Finiteness and homological properties}
\label{sec:transfer}

Throughout this section, \(M\) is a \(w\)-module, except in
the explicitly stated counterexamples.

\begin{lemma}
\label{lem:wft-correspondence}
Let $I \ltimes N$ be a homogeneous ideal of $S = R \ltimes M$.
Then $I \ltimes N$ is of $w$-finite type as an ideal
of $S$ if and only if $I$ is of $w$-finite type as an ideal of $R$ and
$N/IM$ is of $w$-finite type as an $R$-module.
\end{lemma}

\begin{proof}
Suppose first that $I\ltimes N$ is of $w$-finite type.
Choose $(a_1,n_1),\ldots,(a_t,n_t)$ whose $S$-span has the same
localizations as $I\ltimes N$ at every $\mathfrak m\ltimes M$,
where $\mathfrak m\in\wMax(R)$.
Projection shows that the $a_i$ generate $I_{\mathfrak m}$.
The second component of a combination of these pairs is, modulo
$I_{\mathfrak m}M_{\mathfrak m}$, an $R_{\mathfrak m}$-linear
combination of the $n_i$. Their images therefore generate
$(N/IM)_{\mathfrak m}$. Thus $I$ and $N/IM$ are of $w$-finite type.

Conversely, choose $a_1,\ldots,a_t\in I$ generating
$I_{\mathfrak m}$ for every $\mathfrak m\in\wMax(R)$, and choose
$n_1,\ldots,n_u\in N$ whose images generate $(N/IM)_{\mathfrak m}$
at all these primes. The ideal
\[
 B=\sum_i S(a_i,0)+\sum_j S(0,n_j)
\]
is finitely generated, and
\[
 B_{\mathfrak m\ltimes M}
 =I_{\mathfrak m}\ltimes
 \left(\sum_j R_{\mathfrak m}n_j+I_{\mathfrak m}M_{\mathfrak m}\right)
 =(I\ltimes N)_{\mathfrak m\ltimes M}.
\]
Hence $I\ltimes N$ is of $w$-finite type.
\end{proof}

The following is the $w$-version of the classical idealization
criterion \cite[Theorem~4.8]{anderson-winders}. Once the
$w$-operations correspond and $M$ is of $w$-finite type,
the equivalence between $w$-Noetherianity of $R$ and that of $S$
also follows from Xing and Wang's $w$-Eakin--Nagata theorem
\cite[Theorem~2.5]{xing-wang-wnoeth}. We include the direct argument,
which also establishes the necessary finite type condition on $M$.

\begin{theorem}
\label{thm:w-noetherian-transfer}
Let $M$ be a $w$-module. Then
\[
R \ltimes M \text{ is } w\text{-Noetherian}
\iff
\begin{cases}
R \text{ is } w\text{-Noetherian, and}\\
M \text{ is of } w\text{-finite type.}
\end{cases}
\]
\end{theorem}

\begin{proof}
Suppose $S$ is $w$-Noetherian. For every ideal $I$ of $R$,
the ideals $I\ltimes M$ and $0\ltimes M$ are of $w$-finite type
over $S$. Lemma~\ref{lem:wft-correspondence}, applied to these
two ideals, shows that $I$ and $M$ are of $w$-finite type over $R$.
Thus $R$ is $w$-Noetherian.

Conversely, suppose $R$ is $w$-Noetherian and $M$ is of
$w$-finite type. Then $S=R\oplus M$ is a $w$-Noetherian
$R$-module \cite[Theorem~6.8.5(6) and Corollary~6.8.3(1)]{wang-monograph2}.
For any ideal $L$ of $S$, choose a finitely generated
$R$-submodule $B_0\subseteq L$ such that
$(B_0)_{\mathfrak m}=L_{\mathfrak m}$ for all
$\mathfrak m\in\wMax(R)$. Put $B=SB_0\subseteq L$.
This is a finitely generated ideal of $S$, and the canonical
localization isomorphisms give
\[
 B_{\mathfrak m\ltimes M}
 =S_{\mathfrak m\ltimes M}(B_0)_{\mathfrak m}
 =L_{\mathfrak m\ltimes M}
 \qquad(\mathfrak m\in\wMax(R)).
\]
Thus $L$ is of $w$-finite type over $S$, proving the converse.
\end{proof}

\begin{remark}
\label{rem:no-homogeneous-reduction}
For an ideal $L\subseteq R\ltimes M$, put
$I_0=\{r:(r,0)\in L\}$ and $N=\{m:(0,m)\in L\}$.
Its largest homogeneous subideal is $I_0\ltimes N$;
the projection of $L$ to $R$ need not equal $I_0$.
The quotient $L/(I_0\ltimes N)$ need not be GV-torsion.
Take $R=\mathbb Z_{(2)}$, $M=\mathbb F_2$ and
$L=S(2,1)=\{(2a,\bar a):a\in R\}$.
Then the projection is $(2)$, but $I_0=(4)$ and $N=0$.
Thus $I_0\ltimes N=(4)\ltimes0$ and
$L/((4)\ltimes0)\cong\mathbb F_2$.
The local ring $S$ is DW: its maximal ideal is annihilated
by $(0,1)$, so no proper ideal can be a GV-ideal.
The nonzero quotient is therefore not GV-torsion.
\end{remark}

\begin{corollary}
\label{cor:w-noetherian-sufficient}
If $R$ is $w$-Noetherian and $M$ is a $w$-module of $w$-finite type, then
$R \ltimes M$ is $w$-Noetherian.
\end{corollary}

\begin{corollary}
\label{cor:domain-case}
Let $D$ be an integral domain with quotient field $K$ and let $E$ be a
$K$-vector space. The ring $D \ltimes E$ is $w$-Noetherian if and only if $D$
is $w$-Noetherian and $E$ is of $w$-finite type. In particular,
\begin{enumerate}[label=(\arabic*),leftmargin=*,nosep]
\item if $D$ is not a field, then $D\ltimes E$ is $w$-Noetherian
if and only if $D$ is $w$-Noetherian and $E=0$;
\item if $D$ is a field, then $D\ltimes E$ is $w$-Noetherian
if and only if $E$ is finite-dimensional over $D$. In this case
$D\ltimes E$ is Artinian.
\end{enumerate}
\end{corollary}

\begin{proof}
Every quotient field vector space $E$ is a $w$-module and
$E_{\mathfrak p}=E$ for every prime $\mathfrak p$ of $D$.
The equivalence follows from Theorem~\ref{thm:w-noetherian-transfer}.
If $D$ is not a field, choose a nonzero $w$-maximal ideal
$\mathfrak m$ and $0\ne a\in\mathfrak m$.
If $E$ is of $w$-finite type, it is finitely generated over
$D_{\mathfrak m}$. Since $aE=E$, Nakayama's lemma gives $E=0$.
If $D$ is a field, $w$-finite type is finite dimension, and
$D\ltimes E$ is then Artinian.
\end{proof}

\begin{proposition}
\label{prop:dw-transfer}
Let $M$ be a $w$-module. Then $R \ltimes M$ is a DW-ring if and only if
$R$ is a DW-ring.
\end{proposition}

\begin{proof}
The usual maximal-ideal correspondence and
Theorem~\ref{thm:wmax-correspondence} give
\[
 \wMax(S)=\Max(S)\iff\wMax(R)=\Max(R).
\]
This is the DW condition.
\end{proof}

\begin{lemma}
\label{lem:finite-type-restriction}
Let $M$ be a $w$-module and $S=R\ltimes M$.
\begin{enumerate}[label=(\arabic*),leftmargin=*,nosep]
\item A sequence of $S$-modules is $w$-exact over $S$ if and only if
its underlying sequence is $w$-exact over $R$.
\item Let $X$ be an $S$-module. If $M$ is of $w$-finite type over $R$
or $(0\ltimes M)X=0$, then $X$ is of $w$-finite type over $S$ if and
only if its underlying $R$-module is of $w$-finite type.
\end{enumerate}
\end{lemma}

\begin{proof}
By Theorem~\ref{thm:wmax-correspondence} and
\eqref{eq:localization-formula}, localization of an $S$-module at
$\mathfrak m\ltimes M$ is its localization as an $R$-module at
$\mathfrak m$. This proves the assertion about $w$-exactness.
For the corresponding restriction of GV-torsion under the induced
$w$-operation, see \cite[Proposition~2.3(1)]{lin-wang-wphi}.
If a finitely generated $R$-submodule $B\subseteq X$ satisfies
$B_{\mathfrak m}=X_{\mathfrak m}$ for every $\mathfrak m\in\wMax(R)$,
then $SB$ is finitely generated over $S$ and
$(SB)_{\mathfrak m\ltimes M}=X_{\mathfrak m\ltimes M}$ for every such
$\mathfrak m$.

Conversely, choose $x_1,\ldots,x_n\in X$ generating $X$ over $S$
at every $\mathfrak m\ltimes M$. If $M$ is of $w$-finite type,
choose $m_1,\ldots,m_t\in M$ generating $M_{\mathfrak m}$ over
$R_{\mathfrak m}$ for every $\mathfrak m\in\wMax(R)$.
The $R$-submodule generated by
\[
x_i\quad\text{and}\quad (0,m_j)x_i
\qquad(1\le i\le n,\ 1\le j\le t)
\]
has localization $X_{\mathfrak m}$ at every such $\mathfrak m$.
Indeed, for any $z\in X_{\mathfrak m}$, write
$z=\sum_i(a_i,u_i)(x_i/1)$ with $a_i\in R_{\mathfrak m}$
and $u_i\in M_{\mathfrak m}$. Since
$u_i=\sum_j b_{ij}(m_j/1)$ for some $b_{ij}\in R_{\mathfrak m}$,
we have
\[
 z=\sum_i a_i\frac{x_i}{1}
   +\sum_{i,j}b_{ij}\frac{(0,m_j)x_i}{1}.
\]
When $(0\ltimes M)X=0$, the $S$-action on $X$ factors through
$S/(0\ltimes M)\cong R$. Hence, for every
$\mathfrak m\in\wMax(R)$, the same elements $x_1,\ldots,x_n$ also
generate $X_{\mathfrak m}$ as an $R_{\mathfrak m}$-module.
\end{proof}

This gives the analogue for modules of $w$-finite type of the classical
coherence criterion for idealizations by finitely generated modules recalled in
Section~\ref{sec:prelim}; see \cite[p.~55 and Section~4.4]{glaz}.

\begin{theorem}
\label{thm:w-coherent-transfer}
Let $M$ be a $w$-module of $w$-finite type. Then
\[
R\ltimes M\text{ is }w\text{-coherent}
\iff
\begin{cases}
R\text{ is }w\text{-coherent},\\
M\text{ is of }w\text{-finitely presented type over }R.
\end{cases}
\]
\end{theorem}

\begin{proof}
Put $J=0\ltimes M$.
Suppose first that $R$ is $w$-coherent and $M$ is of
$w$-finitely presented type. As an $R$-module,
$S=R\oplus M$ is of $w$-finitely presented type, hence is
$w$-coherent by Proposition~\ref{prop:finite-type-tools}(4).
Let $L$ be a finitely generated ideal of $S$.
Lemma~\ref{lem:finite-type-restriction} makes $L$ a submodule of
$w$-finite type of the $w$-coherent $R$-module $S$.
Thus $L$ is of $w$-finitely presented type over $R$.
Choose an exact sequence
\[
0\longrightarrow K\longrightarrow S^n\longrightarrow L
\longrightarrow0.
\]
The underlying $R$-module $S^n$ is of $w$-finite type, so
Proposition~\ref{prop:finite-type-tools}(2) gives that $K$ is
of $w$-finite type over $R$, and hence over $S$ by
Lemma~\ref{lem:finite-type-restriction}.
Part~(1) of Proposition~\ref{prop:finite-type-tools}, now over
$S$, makes $L$ of $w$-finitely presented type over $S$.
Since this holds for every finitely generated ideal, $S$ is
$w$-coherent.

Conversely, suppose that $S$ is $w$-coherent.
The ideal $J$ is of $w$-finite type over $S$ because $M$ is of
$w$-finite type over $R$. Therefore $J$ is of
$w$-finitely presented type over $S$, and
$S/J\cong R$ is a $w$-coherent $S$-module by
Proposition~\ref{prop:finite-type-tools}(3).
Let $I$ be a finitely generated ideal of $R$, regarded as an
$S$-submodule of $R=S/J$. It is of $w$-finitely presented type
over $S$. In an exact sequence
\[
0\longrightarrow K\longrightarrow R^n\longrightarrow I
\longrightarrow0,
\]
the middle term is finitely generated over $S$.
Proposition~\ref{prop:finite-type-tools}(2) gives that $K$ is
of $w$-finite type over $S$. Since $JK=0$, it is of
$w$-finite type over $R$ by
Lemma~\ref{lem:finite-type-restriction}. Part~(1) of
Proposition~\ref{prop:finite-type-tools} now gives that $I$ is
of $w$-finitely presented type over $R$. Hence $R$ is
$w$-coherent.

Finally, take a $w$-epimorphism $R^n\to M$ and let $H$ be its
kernel. Give $R^n$ and $M$ the $S$-actions through $S/J$.
The sequence
\[
0\longrightarrow H\longrightarrow R^n\longrightarrow M
\longrightarrow0
\]
is $w$-exact over $S$ by
Lemma~\ref{lem:finite-type-restriction}. Its last module is
isomorphic to $J$ as an $S$-module, and so is of
$w$-finitely presented type. Part~(2) of
Proposition~\ref{prop:finite-type-tools} makes $H$ of
$w$-finite type over $S$, hence over $R$ because $JH=0$.
Part~(1) then makes $M$ of $w$-finitely presented type over $R$.
\end{proof}

This is a $w$-version of coherence descent to ring retracts
\cite[Theorem~4.1.5]{glaz}. The proof uses the projection between relation modules and
Lemma~\ref{lem:finite-type-restriction}.

\begin{theorem}\label{thm:coherence-descent}
If \(S\) is \(w\)-coherent, then \(R\) is \(w\)-coherent.
\end{theorem}
\begin{proof}
Let \(I=(a_1,\ldots,a_n)\subseteq R\), and let
\(H=\ker(R^n\to I)\), with the indicated row of generators.
The ideal \(IS\) is finitely generated over \(S\), hence of
\(w\)-finitely presented type. Therefore
\[
 K=\ker\bigl(S^n\longrightarrow IS\bigr)
\]
is of \(w\)-finite type over \(S\).
First components give an \(S\)-linear surjection \(K\to H\),
where \(H\) is acted on through \(S\to R\): every \(h\in H\)
is the first component of \((h,0)\in K\).
Consequently \(H\) is of \(w\)-finite type over \(S\), and hence
over \(R\), since \(JH=0\). The relation criterion over \(R\)
makes \(I\) of \(w\)-finitely presented type. This holds for every
finitely generated ideal \(I\).
\end{proof}

For a $w$-linked homomorphism $\varphi:A\to B$,
\cite[Proposition~2.13]{lin-wang-wphi} asserts that a $B$-module
is $w_\varphi$-finitely generated (respectively, $w_\varphi$-finitely
presented) if and only if its underlying $A$-module is of
$w$-finite type (respectively, of $w$-finitely presented type),
where $w_\varphi$ is the operation induced by $\varphi$.
For a ring homomorphism
$A\to B$, restriction of scalars from $B$ to $A$ need not preserve
finite generation. For example, if $k$ is a field, then $k[t]$ is
free of rank one over itself but infinite-dimensional over $k$.
For $k\to k[t]$ the induced operation is the identity, so
restriction of scalars need not preserve
$w$-finite type.

\begin{theorem}\label{thm:module-fp-comparison}
If \(M\) is of \(w\)-finitely presented type over \(R\), then
an arbitrary \(S\)-module \(X\) is of \(w\)-finitely presented
type over \(S\) if and only if its underlying \(R\)-module is
of \(w\)-finitely presented type.
\end{theorem}
\begin{proof}
The underlying \(R\)-module \(S=R\oplus M\) is of
\(w\)-finitely presented type. Also, \(w\)-finite type over the
two rings is equivalent when \(M\) is of \(w\)-finite type.
Choose an \(S\)-linear \(w\)-epimorphism \(S^n\to X\),
and put \(K=\ker(S^n\to X)\).
If \(X\) is of \(w\)-finitely presented type over \(S\), then
\(K\) is of \(w\)-finite type over \(S\), hence over \(R\).
The \(w\)-exact sequence
\(0\to K\to S^n\to X\to0\), now over \(R\), makes \(X\)
of \(w\)-finitely presented type over \(R\).

Conversely, if $X$ is of $w$-finitely presented type over $R$, it is of
\(w\)-finite type over \(S\), so the chosen map exists.
Since \(S^n\) is of \(w\)-finite type over \(R\), the relation
criterion makes \(K\) of \(w\)-finite type over \(R\), hence
over \(S\). It follows that \(X\) is of \(w\)-finitely presented
type over \(S\).
\end{proof}

A ring $A$ is $w$-coherent if and only if, for every $n,t\ge1$
and every $A$-linear map $\varphi:A^n\to A^t$, the $A$-module
$\ker\varphi$ is of $w$-finite type
\cite[Section~6.11.1]{wang-monograph2}.

\begin{theorem}\label{thm:faithful-coherence}
Suppose there exist \(m_1,\ldots,m_t\in M\), with \(t\ge1\), such that
\[
 \bigcap_{i=1}^t\Ann_R(m_i)=0.
\]
Then \(S\) is \(w\)-coherent if and only if \(R\) is
\(w\)-coherent and \(M\) is of \(w\)-finitely presented type
over \(R\). In particular, this holds if \(\Ann_R(m)=0\) for some \(m\in M\).
\end{theorem}
\begin{proof}
Suppose \(S\) is \(w\)-coherent. The map
\[
 S\longrightarrow S^t,\qquad
 s\longmapsto\bigl(s(0,m_1),\ldots,s(0,m_t)\bigr)
\]
has kernel \(J=0\ltimes M\). Its image is a finitely generated
submodule of \(S^t\), hence of \(w\)-finitely presented type.
Proposition~\ref{prop:finite-type-tools}(2) makes \(J\) of
\(w\)-finite type over \(S\). Since \(J^2=0\),
Lemma~\ref{lem:finite-type-restriction} makes \(M\) of
\(w\)-finite type over \(R\). The result now follows from
Theorem~\ref{thm:w-coherent-transfer}, which also gives the
reverse implication.
\end{proof}

The next lemma applies the standard nilpotent form of Nakayama's
lemma \cite[Theorem~1.3.21(3) and Corollary~1.3.22]{wang-monograph2}
at every $w$-maximal localization. In the proof, finitely many
elements $x_1,\ldots,x_n\in X$ are chosen so that their images
generate $(X/JX)_{\mathfrak m}$ over $R_{\mathfrak m}$ for every
$\mathfrak m\in\wMax(R)$. The ordinary idealization version
is recorded in the discussion of $n=0$ after
\cite[Corollary~4.4.3]{glaz}.

\begin{lemma}
\label{lem:w-patching}
Let $M$ be a $w$-module and set
$S=R\ltimes M$ and $J=0\ltimes M$.
An $S$-module $X$ is of $w$-finite type over $S$ if and only
if $X/JX$ is of $w$-finite type over $R$.
\end{lemma}

\begin{proof}
For modules killed by $J$, finite type over $S$ and over $R$
are equivalent: their actions factor through $S/J=R$ and
their localizations correspond by
Theorem~\ref{thm:wmax-correspondence}.
If $B\subseteq X$ is finitely generated over $S$ and $X/B$ is
GV-torsion, then $(B+JX)/JX$ is finitely generated over $R$ and
$X/(B+JX)$ is GV-torsion. Thus $X/JX$ is of $w$-finite type.

Conversely, choose $\bar x_1,\ldots,\bar x_n\in X/JX$ such that
$(X/JX)/\sum_iR\bar x_i$ is GV-torsion. Lift them to
$x_1,\ldots,x_n\in X$, and put $B=\sum_i Sx_i$ and $W=X/B$.
Then $W/JW$ is GV-torsion over $R$ and, since it is killed
by $J$, over $S$. At every $\mathfrak q=\mathfrak m\ltimes M$
in $\wMax(S)$ this gives
$W_{\mathfrak q}=J_{\mathfrak q}W_{\mathfrak q}$.
Applying $J_{\mathfrak q}$ again and using $J^2=0$ gives
$W_{\mathfrak q}=0$. Hence $W$ is GV-torsion over $S$, and
$X$ is of $w$-finite type because $B$ is finitely generated.
\end{proof}

\begin{theorem}\label{thm:inflated-fp}
Let \(X\) be an arbitrary \(R\)-module, regarded as an
\(S\)-module through \(S\to R\). The following are equivalent:
\begin{enumerate}[label=(\roman*),leftmargin=*]
\item \(X\) is of \(w\)-finitely presented type over \(S\).
\item \(X\) is of \(w\)-finitely presented type over \(R\), and
\(M\otimes_R X\) is of \(w\)-finite type over \(R\).
\end{enumerate}
\end{theorem}
\begin{proof}
Since $JX=0$, Lemma~\ref{lem:w-patching} shows that $X$ is of
$w$-finite type over $S$ if and only if it is of $w$-finite type
over $R$. Thus either (i) or (ii) implies that $X$ is of
$w$-finite type over $R$.
Choose an \(R\)-linear \(w\)-epimorphism \(f:R^n\to X\), put
\(H=\ker f\) and \(B=\operatorname{im}f\), and extend \(f\)
to \(S^n\to X\) through projection. Its kernel, as an
\(R\)-module, is
\[
 K=H\oplus M^n,\qquad
 JK=0\oplus\operatorname{im}(M\otimes_R H\to M^n).
\]
Tensoring \(H\to R^n\to B\to0\) therefore gives
\[
 K/JK\cong H\oplus(M\otimes_R B).
\]
The inclusion \(B\to X\) is a \(w\)-isomorphism.
Its tensor with \(M\) is a \(w\)-isomorphism as well,
as one sees by localizing, where \(B_{\mathfrak m}=X_{\mathfrak m}\).
By Proposition~\ref{prop:finite-type-tools}(1), \(X\) is of \(w\)-finitely presented
type over \(S\) if and only if \(K\) is of \(w\)-finite type over \(S\).
By Lemma~\ref{lem:w-patching} applied to \(K\), this is equivalent
to \(K/JK\) being of \(w\)-finite type over \(R\).
Using $K/JK\cong H\oplus(M\otimes_R B)$ and the
$w$-isomorphism $M\otimes_R B\to M\otimes_R X$, we obtain
\[
 \begin{gathered}
 K/JK\text{ is of }w\text{-finite type over }R
 \\
 \Longleftrightarrow
 \\
 H\text{ and }M\otimes_R X\text{ are both of }w\text{-finite type over }R.
 \end{gathered}
\]
Finally, applying Proposition~\ref{prop:finite-type-tools}(1) to
the $w$-exact sequence
\[
 0\longrightarrow H\longrightarrow R^n\xrightarrow{f}X\longrightarrow0
\]
shows that $H$ is of $w$-finite type over $R$ if and only if
$X$ is of $w$-finitely presented type over $R$.
\end{proof}

\begin{corollary}\label{cor:quotient-fp}
The \(S\)-module \(R=S/J\) is of \(w\)-finitely presented type
if and only if \(M\) is of \(w\)-finite type over \(R\).
If \(I\subseteq R\) is finitely generated, then \(R/I\) is of
\(w\)-finitely presented type over \(S\) if and only if
\(M/IM\) is of \(w\)-finite type over \(R\).
\end{corollary}
\begin{proof}
The $R$-module $R$ is finitely presented, and so is $R/I$ because
$I$ is finitely generated. Using
\[
 M\otimes_R R\cong M,\qquad M\otimes_R(R/I)\cong M/IM,
\]
apply Theorem~\ref{thm:inflated-fp} with $X=R$ and $X=R/I$,
respectively, to obtain the two assertions.
\end{proof}

Combining \cite[Theorem~6.8.5 and Section~6.11.1]{wang-monograph2}
with the transfer theorems yields the following corollary.

\begin{corollary}
\label{cor:noetherian-coherent}
\begin{enumerate}[label=(\arabic*),leftmargin=*,nosep]
\item Let $M$ be a $w$-module and $S=R\ltimes M$. If $S$ is
$w$-Noetherian, then $S$ is $w$-coherent. The converse holds when
$R$ is $w$-Noetherian and $M$ is of $w$-finite type.
\item Over a $w$-Noetherian ring,
the classes of modules of $w$-finitely presented type and of
$w$-finite type coincide.
\end{enumerate}
\end{corollary}

\begin{remark}
\label{rem:chang-kim-wff}
The $w$-FF property in \cite{chang-kim,kim-lim} requires
$w$-flat $w$-ideals to be of $w$-finite type. It is a
different condition from $w$-coherence, which requires every ideal
of $w$-finite type to be of $w$-finitely presented type.
For example, the Chang--Kim transfer theorem gives the
$w$-FF property for $D\ltimes E$ when $D$ is a DVR and
$E$ is a nonzero vector space over its quotient field.
However, for any $0\ne e\in E$, the annihilator
\[
 \Ann_{D\ltimes E}(0,e)=0\ltimes E
\]
is not of $w$-finite type. Hence $D\ltimes E$ is not
$w$-coherent by Example~\ref{ex:domain-vector-spaces}.
\end{remark}

The next theorem and its canonical module corollaries are
$w$-local applications of Lemma~\ref{lem:classical-reiten-fossum}.

\begin{theorem}
\label{thm:w-gorenstein}
Let $M$ be a $w$-module.
Then $R\ltimes M$ is $w$-Gorenstein if and only if both of the
following conditions hold:
\begin{enumerate}[label=(\alph*),leftmargin=*]
\item $R$ is $w$-Noetherian and $M$ is of $w$-finite type;
\item for every $\mathfrak m\in\wMax(R)$, either
  \begin{enumerate}[label=(\roman*),leftmargin=*]
  \item $M_{\mathfrak m}=0$ and $R_{\mathfrak m}$ is Gorenstein, or
  \item $M_{\mathfrak m}\ne0$, $R_{\mathfrak m}$ is
  Cohen--Macaulay, and $M_{\mathfrak m}\cong\omega_{R_{\mathfrak m}}$.
  \end{enumerate}
\end{enumerate}
\end{theorem}

\begin{proof}
By Theorem~\ref{thm:w-noetherian-transfer}, the global
$w$-Noetherian condition on $S$ is equivalent to $R$ being
$w$-Noetherian and $M$ being of $w$-finite type. Under (a),
$A=R_{\mathfrak m}$ is Noetherian and
$E=M_{\mathfrak m}$ is finitely generated for every
$\mathfrak m\in\wMax(R)$.
By Theorem~\ref{thm:wmax-correspondence} and
\eqref{eq:localization-formula}, it remains to require that
$R_{\mathfrak m}\ltimes M_{\mathfrak m}$ be Gorenstein for every
$\mathfrak m\in\wMax(R)$.
If $E=0$, this is exactly alternative (i) in condition (b).
If $E\ne0$, Lemma~\ref{lem:classical-reiten-fossum} makes this
equivalent to alternative (ii) in condition (b).
\end{proof}

\begin{remark}
\label{rem:w-gorenstein-global}
Even when $R\ltimes M$ is Gorenstein and $M\ne0$, one can have
$M_{\mathfrak m}=0$ for some $\mathfrak m\in\wMax(R)$.
For example, let $k$ be a field, $R=k\times k$, and $M=k\times0$.
Then $R\ltimes M\cong k[\varepsilon]/(\varepsilon^2)\times k$ is
Gorenstein. At $\mathfrak m=k\times0$, we have
$M_{\mathfrak m}=0$ and $R_{\mathfrak m}\cong k$, as in
Theorem~\ref{thm:w-gorenstein}(b)(i).
Both rings are zero-dimensional DW-rings, so the $w$-module
hypothesis holds.

The global finiteness hypotheses also cannot be replaced by
the local canonical module condition. Let $R=\prod_{i\ge1}k$
and $M=R$. The module $M$ is of $w$-finite type, generated by $1$.
Every prime localization of $R$ is a field and
$M_{\mathfrak m}=R_{\mathfrak m}$ is canonical there.
Nevertheless $R$ is not Noetherian. The ideal
$\bigoplus_{i\ge1}k$ is not finitely generated. Since $R$ is
a DW-ring, it is not $w$-Noetherian, and neither is $R\ltimes R$.
Hence $R\ltimes R$ is not $w$-Gorenstein.
\end{remark}

\begin{corollary}
\label{cor:w-fossum}
If $R$ is $w$-Cohen--Macaulay and $C$ is a $w$-canonical module
of $w$-finite type,
then $R\ltimes C$ is $w$-Gorenstein.
\end{corollary}

\begin{proof}
Apply Theorem~\ref{thm:w-gorenstein} with $M=C$, using
$0\ne C_{\mathfrak m}\cong\omega_{R_{\mathfrak m}}$ for every
$\mathfrak m\in\wMax(R)$.
\end{proof}

\begin{corollary}[Classical Reiten--Fossum theorem {\cite{fossum,reiten}}]
\label{cor:classical-fossum}
If $A$ is a Noetherian Cohen--Macaulay local ring with canonical module
$\omega_A$, then $A\ltimes\omega_A$ is Gorenstein local.
More generally, if $R$ is Noetherian and has a finitely
generated module $C$ whose localization at every maximal
ideal is canonical, then $R\ltimes C$ is Gorenstein.
\end{corollary}
\begin{proof}
Apply Lemma~\ref{lem:classical-reiten-fossum} at every maximal ideal.
\end{proof}

The following examples distinguish the finiteness hypotheses and
the local canonical module condition.

\begin{example}
\label{ex:classical-setting}
Let $R$ be Noetherian and let $M$ be a finitely generated
$w$-module, for instance $M=R^n$. Then $M$ is finitely presented,
so the transfer theorems give both $w$-Noetherianity and
$w$-coherence of $R\ltimes M$. In fact, $R\ltimes M$ is Noetherian
\cite[Theorem~4.8]{anderson-winders}.
\end{example}

\begin{example}
\label{ex:z-q}
Let $S=\mathbb Z\ltimes\mathbb Q$. Here $\mathbb Q$ is flat
and hence a $w$-module.
The ring $\mathbb Z$ is DW, so $S$ is DW by
Proposition~\ref{prop:dw-transfer}.
Since $p\mathbb Q=\mathbb Q$ for every prime $p$, Nakayama's
lemma shows that $\mathbb Q$ is not finitely generated over
$\mathbb Z_{(p)}$. Thus it is not of $w$-finite type, and $S$
is not $w$-Noetherian.

Also,
\[
 \Ann_S(0,1)=0\ltimes\mathbb Q
\]
is not finitely generated: the $S$-span of finitely many
$(0,q_i)$ is $0\ltimes\sum_i\mathbb Zq_i$.
Thus $S$ is not coherent \cite[Theorem~2.8(1) and Example~2.11]{kabbaj-mahdou-coherent},
and hence not $w$-coherent. Even over $\mathbb Z$, flatness of the
extension module alone guarantees neither $w$-Noetherianity nor
$w$-coherence of the extension.
\end{example}

\begin{remark}
\label{rem:zckz-comparison}
Let $D$ be a domain with quotient field $K$, and write
$K[\theta]=K[t]/(t^2)$, where $\theta$ is the class of $t$.
The map $\pi:K[\theta]\to K$, $a+b\theta\mapsto a$, gives the pullback square
\[
\begin{tikzcd}
R_2=D+\theta K \ar[r,hook] \ar[d,two heads,"\pi|_{R_2}"'] & K[\theta] \ar[d,two heads,"\pi"] \\
D \ar[r,hook] & K
\end{tikzcd}
\]
Thus $R_2=\pi^{-1}(D)\cong D\ltimes K$ via $d+b\theta\mapsto(d,b)$.
The GV-ideal and $w$-spectral formulas in \cite[Lemma~2.1]{zckz}
already cover this case of Section~\ref{sec:internal}.

In contrast, the coherence assertions of \cite[Theorems~1.2 and~2.4]{zckz}
fail for $D=\mathbb Z$ and $K=\mathbb Q$. Indeed,
\[
 \Ann_{R_2}(\theta)=\theta\mathbb Q
 \ne\theta\mathbb Z=\theta R_2.
\]
The ideal $\theta\mathbb Q$ is not finitely generated, so $R_2$ is not
coherent \cite[Theorem~2.8(1)]{kabbaj-mahdou-coherent}.
Since $R_2$ is a DW-ring, it is not $w$-coherent either
(Example~\ref{ex:z-q}). The equality
$\Ann_{R_2}(\theta)=\theta R_2$ used in
\cite[Theorem~2.4, proof on p.~3256]{zckz} is contradicted by this computation.

For $K[\theta]=K[t]/(t^n)$ with $n\ge3$, the ideal $\theta K[\theta]$
is not square-zero, so the idealization results do not directly apply
to $D+\theta K[\theta]\to D$.
\end{remark}

\begin{example}
\label{ex:field-base}
Let $K$ be a field and $E$ a finite-dimensional $K$-vector space.
Then $K\ltimes E$ is always Artinian, but it is Gorenstein if and
only if $\dim_K E\le 1$. Indeed, $K\ltimes E$ is a finite-dimensional
local $K$-algebra, so it is Artinian and hence $w$-Noetherian and
$w$-coherent. If $E=0$, then $K\ltimes E\cong K$. If $E\ne0$,
its socle is $0\ltimes E$, which is one-dimensional precisely when
$\dim_K E=1$. This is the field case of
Lemma~\ref{lem:classical-reiten-fossum}.
\end{example}

\begin{example}
\label{ex:domain-vector-spaces}
Let $D$ be a domain with quotient field $K$, let $E$ be a nonzero
$K$-vector space, and put $S=D\ltimes E$. Then $S$ is $w$-coherent
if and only if $D$ is a field and $\dim_K E<\infty$.

For necessity, we adapt the annihilator argument of
\cite[Theorem~2.8(1)]{kabbaj-mahdou-coherent}; see also
\cite[Proposition~2.5]{adarbeh-kabbaj}. Choose $0\ne e\in E$ and put $I=S(0,e)$.
Since $(d,u)(0,e)=(0,de)$, we have $\Ann_S(I)=0\ltimes E$.
If $S$ is $w$-coherent, then $I$
is of $w$-finitely presented type, so
Proposition~\ref{prop:finite-type-tools}(1) implies
that $0\ltimes E$ is of $w$-finite type over $S$.
Its action factors through $D$; the spectral
correspondence then makes $E$ of $w$-finite type over
$D$. The argument of Corollary~\ref{cor:domain-case}
forces $D$ to be a field and $E$ finite-dimensional.
Conversely, if $D$ is a field and $\dim_D E<\infty$, then $S$ is
finite-dimensional over $D$, hence Artinian and Noetherian. Every
finitely generated ideal of $S$ is therefore finitely presented, so
Proposition~\ref{prop:finite-type-tools}(4) yields that $S$ is
$w$-coherent.
\end{example}

\begin{example}
\label{ex:wm-hyp-exact}
Let $R=k[x,y]$, $\mathfrak m=(x,y)$, and $M=R/\mathfrak m=k$.
Since $\mathfrak m\in\GV(R)$ and $\mathfrak mM=0$ with $M\ne0$,
$M$ is not GV-torsion-free and hence is not a $w$-module.
The next proposition computes $\GV(R\ltimes M)$ and shows that
$\wMax(R\ltimes M)$ and
$\{\mathfrak p\ltimes M:\mathfrak p\in\wMax(R)\}$ are incomparable.
Thus the $w$-module hypothesis in
Theorem~\ref{thm:wmax-correspondence} cannot be omitted.
\end{example}

\begin{proposition}
\label{prop:exactness}
Let $R=k[x,y]$, $\mathfrak m=(x,y)$, $M=R/\mathfrak m=k$,
and $S=R\ltimes M$.
\begin{enumerate}[label=(\arabic*),leftmargin=*,nosep]
\item The GV-ideals of $S$ are given by
\begin{equation}\label{eq:gv-characterization}
\GV(S)=\{I\ltimes k:I\text{ finitely generated},\
(R:I)=R,\ I+\mathfrak m=R\}.
\end{equation}
\item $M$ is not a $w$-module and the sets
$\wMax(S)$ and
$\{\mathfrak p\ltimes k:\mathfrak p\in\wMax(R)\}$
are incomparable under inclusion.
\end{enumerate}
\end{proposition}

\begin{proof}
Since $x,y$ are coprime in the UFD $R$, a fraction $q$ with
$qx,qy\in R$ belongs to $R$. Hence $(R:\mathfrak m)=R$ and
$\mathfrak m\in\GV(R)$. This ideal kills the nonzero module
$M$, so $M$ is not GV-torsion-free.

Let $G\in\GV(S)$ and let $I$ be its projection to $R$.
The ideal $I$ is finitely generated. If $I\subseteq\mathfrak m$,
then $(0\ltimes k)G=0$, contradicting faithfulness of $G$.
Thus $I+\mathfrak m=R$. Choose $(c,u)\in G$ with
$c\equiv1\pmod{\mathfrak m}$. Multiplication by $(0,v)$ shows
that $(0,v)\in G$ for every $v\in k$. Consequently
$G=I\ltimes k$.

Take $q\in(R:I)$ and put $\lambda=\overline{qc}\in k$.
For every $a\in I$, comparison of the images of
$(qa)c=(qc)a$ modulo $\mathfrak m$ gives
$\overline{qa}=\lambda\bar a$.
It follows that
\[
f:I\ltimes k\longrightarrow S,\qquad
f(a,v)=(qa,\lambda v)
\]
is $S$-linear. Indeed, the second components in
$f((r,u)(a,v))=(r,u)f(a,v)$ agree by the displayed identity.
Since $G$ is a GV-ideal, $f$ is multiplication by an element
$(r_0,v_0)\in S$. Comparing first components and using $I\ne0$
gives $q=r_0\in R$. Hence $(R:I)=R$.

Conversely, let $I$ satisfy the conditions on the right of
\eqref{eq:gv-characterization}. Since $Ik=k$, the ideal
$G=I\ltimes k$ is generated by the pairs $(a_i,0)$ for a
finite generating set $a_i$ of $I$.
The ideal $G$ is faithful. Indeed, if $(r,v)\in\Ann_S(G)$,
then $rI=0$, hence $r=0$. Choose $c\in I$ with $\bar c=1$.
The equality $(0,v)(c,0)=(0,v)=0$ then gives $v=0$.
For an $S$-linear map $f:G\to S$, write
$f(a,0)=(\alpha(a),\beta(a))$.
The map $\alpha:I\to R$ is multiplication by some $r_0\in R$
because $(R:I)=R$.
Multiplication by $(0,v)$ and the identity $\bar c=1$ give
$f(0,v)=(0,r_0v)$.
Also $I\cap\mathfrak m=I\mathfrak m$, so
$I/\mathfrak m I\cong k$ and
$\beta(a)=\bar a v_0$ for some $v_0\in k$.
Thus $f$ is multiplication by $(r_0,v_0)$, proving that
$G\in\GV(S)$.

The maximal ideal $\mathfrak m\ltimes k$ contains no GV-ideal
of $S$, so it is $w$-maximal, whereas $\mathfrak m$ is not
a $w$-ideal of $R$. On the other hand, $(x)$ is $w$-maximal
in $R$: it is principal and hence $w$-closed, and adjoining
any $g\notin(x)$ produces the GV-ideal $(x,g)$.
The prime $(x)\ltimes k$ is $w$-closed in $S$ by
\eqref{eq:gv-characterization}, but is properly contained in
the proper $w$-ideal $\mathfrak m\ltimes k$.
It is therefore not $w$-maximal. These two primes prove
incomparability.
\end{proof}

\begin{remark}
\label{rem:mechanism}
More generally, if $0\ne m_0\in M$ and
$\pi(G)\subseteq\Ann_R(m_0)$, then $(0,m_0)G=0$.
Hence $G$ cannot be a GV-ideal.
\end{remark}

\Needspace{14\baselineskip}
\begin{example}\label{ex:local-noetherian-not-global}
Let $R=\mathbb Z$ and
\[
M=\bigoplus_{p\text{ prime}}\mathbb Z/p\mathbb Z.
\]
The ring $R$ is a DW-ring, so every $R$-module is a $w$-module
and $w$-finite type is ordinary finite generation. For each prime
$p$, localization gives $M_{(p)}\cong\mathbb Z/p\mathbb Z$.
The summands indexed by $q\ne p$ vanish because $q$ is invertible
in $\mathbb Z_{(p)}$. Hence
\[
S_{(p)\ltimes M}\cong
\mathbb Z_{(p)}\ltimes(\mathbb Z/p\mathbb Z)
\]
is Noetherian, since it is finitely generated over the Noetherian
ring $\mathbb Z_{(p)}$. Nevertheless $M$ is not finitely generated.
Every finite subset is supported at only finitely many primes.
Theorem~\ref{thm:w-noetherian-transfer} therefore shows that
$S=\mathbb Z\ltimes M$ is not $w$-Noetherian. Thus Noetherianity
of every $w$-maximal localization of the extension, even over a
Noetherian base, does not replace the global $w$-finite-type
condition on the extension module.
\end{example}

\begin{example}\label{ex:localization-obstruction}
Let \(R=k[x,y]\), \(M=R[1/x]\), and \(S=R\ltimes M\).
The module \(M\) is flat, hence a \(w\)-module, and the
\(w\)-operation on \(R\) is nontrivial since \((x,y)\in\GV(R)\).
The prime \((x)\) is \(w\)-maximal. It is principal and
\(w\)-closed, whereas adjoining any \(g\notin(x)\) produces
the GV-ideal \((x,g)\). At this prime,
\[
 M_{(x)}=\operatorname{Frac}(R)
\]
is not finitely generated over the DVR \(R_{(x)}\).
For example, valuations of finitely many prospective generators
are bounded below, while the powers \(x^{-n}\) have unbounded
negative valuations. Thus \(M\) is not of \(w\)-finite type.

The \(S\)-module \(R\) is not of \(w\)-finitely presented type
by Corollary~\ref{cor:quotient-fp}.
In contrast, for every \(n\ge1\),
\[
 X_n=R/(x^n),\qquad M\otimes_R X_n=0,
\]
so \(X_n\) is of \(w\)-finitely presented type over \(S\).

For \(X=S\), we have \(X/JX\cong R\) and
\(\ker\mu_S=0\) (equivalently, \(\Tor_1^S(R,S)=0\)).
Since \(R\) is free of rank one over itself,
Theorem~\ref{thm:w-tor-lifting}(2) shows that \(S\) is of
\(w\)-finitely presented type over \(S\), although its underlying
\(R\)-module is not of \(w\)-finite type. Since \(1\in M\) has zero annihilator,
Theorem~\ref{thm:faithful-coherence} implies that \(S\) is not
\(w\)-coherent. For this fixed $M$, the tensor products satisfy
\[
M\otimes_R R\cong M,\qquad M\otimes_R R/(x^n)=0\quad(n\ge1).
\]
Since $M$ is not of $w$-finite type, Theorem~\ref{thm:inflated-fp}
shows that $R$ is not of
$w$-finitely presented type over $S$, whereas every $R/(x^n)$ is.
\end{example}

\begin{example}\label{ex:divisible-torsion}
Let \(R=\mathbb Z\), \(M=\mathbb Q/\mathbb Z\), and
\(S=\mathbb Z\ltimes(\mathbb Q/\mathbb Z)\).
Every nonzero ideal of \(\mathbb Z\) is principal, and
\(n\mathbb Z\) is a GV-ideal precisely when \(n=1\), up to sign.
Thus \(R\) is a DW-ring and every \(R\)-module is a
\(w\)-module. The characterization of GV-ideals then shows that \(S\)
is a DW-ring as well.

The coherence of \(S\) follows from \cite[Proposition~2.5]{adarbeh-kabbaj}.
Let \(G\) be a finitely generated ideal of \(S\). If its projection is \(n\mathbb Z\)
with \(n>0\), then \(G\) contains \((n,u)\) for some \(u\).
Divisibility, \(nM=M\), gives \(0\ltimes M\subseteq G\).
Consequently
\[
 G=n\mathbb Z\ltimes M=S(n,0).
\]
If its projection is zero, \(G=0\ltimes N\) for a finite
subgroup \(N\subseteq\mathbb Q/\mathbb Z\); such a subgroup
is cyclic, so this ideal is principal too.
Intersections of two ideals of nonzero projection are
\(\operatorname{lcm}(n,m)\mathbb Z\ltimes M\).
The intersection of one such ideal with \(0\ltimes N\) is
\(0\ltimes N\); intersections of two finite subgroups are finite.
All these intersections are finitely generated.

If \(a\ne0\), direct multiplication gives
\[
 \Ann_S(a,u)=0\ltimes M[a],
\]
where \(M[a]\) is cyclic of order \(|a|\).
If \(u\) has order \(n>0\), then
\[
 \Ann_S(0,u)=n\mathbb Z\ltimes M=S(n,0).
\]
The zero element has annihilator \(S\).
Thus $\Ann_S(s)$ is finitely generated for every $s\in S$,
and $G\cap H$ is finitely generated for any two finitely generated
ideals $G,H\subseteq S$. These two conditions characterize
coherent rings, so $S$ is coherent.

Nevertheless \(J=0\ltimes M\) is not finitely generated.
A finite subgroup cannot contain elements of every finite order.
Hence \(S\) is coherent and \(w\)-coherent, but is neither
Noetherian nor \(w\)-Noetherian. Moreover \(R=S/J\) is not
of \(w\)-finitely presented type as an \(S\)-module, while
each \(R/nR\), \(n\ge2\), is finitely presented over \(S\),
because \(M/nM=0\).
The module \(M\) is faithful, but every finite subset has a common
nonzero integer annihilator. Hence no finite family satisfies the
hypothesis of Theorem~\ref{thm:faithful-coherence}. This example
shows that faithfulness of \(M\) alone does not imply that \(M\) is of
\(w\)-finitely presented type over \(R\), even when \(S\) is coherent.
\end{example}

\begin{example}\label{ex:tor-converse-fails}
Let \(k\) be a field and \(V\) an infinite-dimensional
\(k\)-vector space. Set \(R=k\ltimes V\), choose
\(0\ne a\in V\subseteq R\), and let
\(S=R\ltimes R\), with square-zero ideal \(J=0\ltimes R\).
Both rings are zero-dimensional, so are DW-rings;
in particular \(M=R\) is a \(w\)-module.
The base \(R\) is not coherent \cite[Theorem~2.6(2)]{kabbaj-mahdou-coherent}, since
\(\Ann_R(a)=V\) is not finitely generated.

Put \(G=S(a,1)\) and \(X=S/G\). The module \(X\) has an
ordinary finite presentation, and
\(X/JX\cong R/aR\) does too. However,
\[
 G\cap J=0\ltimes V,\qquad JG=0\ltimes aR,
 \qquad
 \Tor_1^S(R,X)\cong (G\cap J)/JG\cong V/ka.
\]
For the first equality, \((r,u)(a,1)=(ra,r+au)\) lies in
\(J\) precisely when \(r\in V\), and
\[
 \{r+au:r\in V,\ u\in R\}=V+aR=V.
\]
The second equality follows by multiplying by \(J\). Tensoring
\(0\to G\to S\to X\to0\) over \(S\) with \(R=S/J\) gives
\[
 \begin{gathered}
 0\longrightarrow\Tor_1^S(R,X)\longrightarrow
 R\otimes_S G\cong G/JG\xrightarrow{\bar\iota}R
 \longrightarrow X/JX\longrightarrow0,\\
 \bar\iota(g+JG)=g+J.
 \end{gathered}
\]
Hence \(\Tor_1^S(R,X)\cong\ker\bar\iota=(G\cap J)/JG\).
Since \(V/ka\) is infinite-dimensional and is annihilated by
\(V\), it is not finitely generated over \(R\).
Thus \(X\) is finitely presented over \(S\), but \(\Tor_1^S(R,X)\)
is not of \(w\)-finite type over \(R\). Hence the \(w\)-coherence assumption
in Theorem~\ref{thm:w-tor-lifting}(3) cannot be dropped.

Over the same ring \(R=k\ltimes V\), set \(E=R/V=k\) and
\(T=R\ltimes E\). Here \(E\) is a cyclic \(w\)-module,
but it is not finitely presented, since the kernel \(V\) of
\(R\to k\) is not finitely generated.
The module \(T\) is free of rank one over itself, while its
underlying \(R\)-module \(R\oplus k\) is not finitely presented:
otherwise its direct summand \(k\) would be finitely presented.
Both rings are DW-rings. Thus the hypothesis in
Theorem~\ref{thm:module-fp-comparison} that the extension module
be of \(w\)-finitely presented type cannot be weakened to
\(w\)-finite type. Indeed, \(T\) is of \(w\)-finitely presented type over
itself, whereas its underlying \(R\)-module is not.
\end{example}

\section{Square-zero cartesian squares and lifting}
\label{sec:patching}

The split square-zero construction is classical
\cite[Section~1]{fgr}; see also \cite[Section~4.4]{glaz}.
We distinguish its common kernel from the conductor.
The $w$-statements below assume that $M$ is a $w$-module over $R$.

\begin{theorem}
\label{thm:square-zero-square}
Let $R\subseteq F$ be an extension of rings and let $M$ be an
$F$-module. Put $S=R\ltimes M$, $T=F\ltimes M$, and
$J=0\ltimes M$.
\begin{enumerate}[label=(\arabic*),leftmargin=*,nosep]
\item The square
\[
\begin{tikzcd}
S \ar[r,hook] \ar[d,two heads] & T \ar[d,two heads]\\
R \ar[r,hook] & F
\end{tikzcd}
\]
is cartesian with common kernel $J$ and $J^2=0$.
\item The full conductor is
\[
(S:T)=(R:F)\ltimes M,
\qquad (R:F)=\{r\in R:rF\subseteq R\}.
\]
In particular the common kernel equals the conductor if and
only if $(R:F)=0$.
\item If $F=U^{-1}R$ is an injective localization, then
$T\cong U^{-1}S\cong S\otimes_R F$ and $T$ is flat over $S$.
\item Conversely, let $q:A\twoheadrightarrow F$ be a surjective ring
homomorphism, put $I=\ker q$ and $B=q^{-1}(R)$, and suppose $I^2=0$.
If there is a ring homomorphism $\sigma:F\to A$ with
$q\circ\sigma=\operatorname{id}_F$, give $I$ the $F$-module structure
$f\cdot i=\sigma(f)i$. Then the ring isomorphism
\[
F\ltimes I\xrightarrow{\ \sim\ }A,\qquad(f,i)\longmapsto\sigma(f)+i,
\]
restricts to an isomorphism $R\ltimes I\cong B$.
Under these isomorphisms, the square in \textup{(1)} with $M=I$ is
identified with the following cartesian square:
\[
\begin{tikzcd}
B \ar[r,hook] \ar[d,two heads,"q|_B"'] & A \ar[d,two heads,"q"]\\
R \ar[r,hook] & F.
\end{tikzcd}
\]
\end{enumerate}
\end{theorem}

\begin{proof}
\textup{(1)} An element of $R\times_F T$ has the form $(r,(r,m))$, giving the
isomorphism with $S$. The kernel and the square-zero property
follow from the multiplication in $T$.

\textup{(2)} For $(r,m)\in S$, all products with $(f,n)\in T$ lie in $S$
exactly when $rf\in R$ for every $f\in F$; the second component
always lies in $M$. This gives the conductor formula.

\textup{(3)} If $F=U^{-1}R$, the $F$-module structure on $M$ identifies
$U^{-1}M$ with $M$, so $U^{-1}S=F\ltimes M$.
The tensor-product identification and flatness follow from
localization.

\textup{(4)} Since $I^2=0$, the displayed map preserves multiplication.
Every $a\in A$ has the unique expression
\[
a=\sigma(q(a))+\bigl(a-\sigma(q(a))\bigr),\qquad
 a-\sigma(q(a))\in I.
\]
Thus the map is a ring isomorphism. Moreover,
$q(\sigma(f)+i)=f$, so it sends $R\ltimes I$ onto $q^{-1}(R)=B$
and commutes with the quotient maps to $R$ and $F$.
\end{proof}

\begin{remark}
\label{rem:square-zero-wlink-fails}
The theorem does not assert a bijection between the $w$-maximal
spectra of $R$ and $F$. For the flat extension
$\mathbb Z\subset\mathbb Z[1/2]$, no prime ideal of $\mathbb Z[1/2]$
contracts to $2\mathbb Z$, since $2$ is a unit in $\mathbb Z[1/2]$.
Nor does an arbitrary $F$-module satisfy
$M\otimes_R F\cong M$; for $R=k$, $F=k[t]$ and $M=F$, the
multiplication map $F\otimes_k F\to F$ has the nonzero kernel
element $t\otimes1-1\otimes t$. Thus the localization assumption
in Theorem~\ref{thm:square-zero-square}(3) cannot be omitted.
\end{remark}

\begin{remark}
\label{rem:range}
The condition that $F/R$ is GV-torsion does not imply $F=R$.
Indeed, for a proper GV-ideal $I$ of $R$, let $F=R\ltimes R/I$
and embed $R$ in $F$ by $r\mapsto(r,0)$. Then $F/R\cong R/I\ne 0$
as $R$-modules and $I(F/R)=0$.
For any extension $R\subseteq F$ with $F/R$ GV-torsion, if $F$ is
GV-torsion-free as an $R$-module, then $F=R$.
Indeed, for $f\in F$ choose $I\in\GV(R)$ with $If\subseteq R$.
The map $a\mapsto af$ from $I$ to $R$ is multiplication by
some $r\in R$, so $I(f-r)=0$ and $f=r$.
\end{remark}

Let $M$ be a $w$-module over $R$, set $S=R\ltimes M$ and $J=0\ltimes M$, and let $X$ be an $S$-module. Put $Q=X/JX$ and define the multiplication map
\[
\mu_X:M\otimes_R Q\longrightarrow JX,
\qquad m\otimes\bar x\longmapsto (0,m)x.
\]
It is well-defined because $J^2=0$, and it is surjective by the
definition of $JX$.

\begin{lemma}[{\cite[Corollary~4.4.3]{glaz}}]
\label{lem:tor-multiplication}
There is a natural isomorphism
\[
\Tor_1^S(R,X)\cong\ker\mu_X.
\]
\end{lemma}

\begin{proof}
Tensor $0\to J\to S\to R\to0$ with $X$ over $S$.
The resulting exact sequence identifies $\Tor_1^S(R,X)$ with
the kernel of $J\otimes_S X\to X$. Since $J^2=0$,
$J\otimes_S X\cong M\otimes_R X/JX$, and this map is $\mu_X$
followed by the inclusion $JX\hookrightarrow X$.
\end{proof}

Roos's criterion for finite presentation is recorded in
\cite[Theorem~4.4.2 and Corollary~4.4.3]{glaz}. The following
theorem gives a criterion for $w$-finite presentation in terms of
$X/JX$ and $\ker\mu_X$; the converse for $\ker\mu_X$ assumes
that $R$ is $w$-coherent.

\begin{theorem}\label{thm:w-tor-lifting}
Let \(M\) be a \(w\)-module, \(S=R\ltimes M\), \(J=0\ltimes M\),
and \(X\) an arbitrary \(S\)-module.
\begin{enumerate}[label=(\arabic*),leftmargin=*,nosep]
\item If \(X\) is of
\(w\)-finitely presented type over \(S\), then \(X/JX\) is
of \(w\)-finitely presented type over \(R\).
\item If \(X/JX\) is of
\(w\)-finitely presented type over \(R\) and
\(\ker\mu_X\) is of \(w\)-finite type over \(R\), then
\(X\) is of \(w\)-finitely presented type over \(S\).
\item If \(R\) is \(w\)-coherent and \(X\) is of \(w\)-finitely
presented type over \(S\), then \(\ker\mu_X\) is of \(w\)-finite
type over \(R\).
\end{enumerate}
\end{theorem}
\begin{proof}
\textup{(1)} Let
\[
0\longrightarrow K\longrightarrow S^n\longrightarrow X\longrightarrow0
\]
be a finite \(w\)-presentation of \(X\). Tensoring this sequence with \(R\)
and checking right exactness at the corresponding localizations gives a
\(w\)-exact sequence
\[
K/JK\longrightarrow R^n\longrightarrow X/JX\longrightarrow0.
\]
Put
\[
L=\ker\bigl(R^n\longrightarrow X/JX\bigr).
\]
Then \(L\) is a quotient of \(K/JK\), so \(L\) is of \(w\)-finite type over
\(R\). Since \(R^n\) is finitely generated free, Proposition~\ref{prop:finite-type-tools}(1)
implies that \(X/JX\) is of \(w\)-finitely presented type over \(R\).

\textup{(2)} Lemma~\ref{lem:w-patching} gives an \(S\)-linear
\(w\)-epimorphism \(f:S^n\to X\). Put \(B=\operatorname{im}f\),
\(K=\ker f\), and \(Q=X/JX\).
The inclusion \(B\to X\) is an isomorphism at all corresponding
localizations, so it induces \(w\)-isomorphisms on tensor
products with \(R\) and on \(\Tor_1^S(R,-)\).
Tensor the exact sequence \(0\to K\to S^n\to B\to0\)
with \(R\). It follows that
\[
 0\longrightarrow T\longrightarrow K/JK
 \longrightarrow R^n\longrightarrow Q\longrightarrow0
 \tag{*}
\]
is \(w\)-exact, where \(T=\Tor_1^S(R,B)\) is \(w\)-isomorphic
to \(\Tor_1^S(R,X)\). Put $L=\ker(R^n\to Q)$.
Since $Q$ is of $w$-finitely presented type over $R$,
Proposition~\ref{prop:finite-type-tools}(1), applied to the $w$-exact sequence
\[
0\longrightarrow L\longrightarrow R^n\longrightarrow Q\longrightarrow0,
\]
shows that $L$ is of $w$-finite type over $R$.
The \(w\)-exact sequence
\(0\to T\to K/JK\to L\to0\) and finite type of \(T,L\)
make $K/JK$ of $w$-finite type. Lemma~\ref{lem:w-patching}
then gives that $K$ is of $w$-finite type over $S$.
Applying Proposition~\ref{prop:finite-type-tools}(1) to the $w$-exact sequence
$0\to K\to S^n\xrightarrow{f}X\to0$ shows that $X$ is of
$w$-finitely presented type over $S$.

\textup{(3)} Suppose \(R\) is \(w\)-coherent and \(X\) is of
\(w\)-finitely presented type over \(S\).
By \textup{(1)}, \(Q\) is of \(w\)-finitely presented
type over \(R\). Choose \(f\) as above. Its kernel \(K\) is of \(w\)-finite type
over \(S\), so \(K/JK\) is of \(w\)-finite type over \(R\).
The module \(L\) in (*) is of \(w\)-finite type, and is a
submodule of the finite free \(R\)-module \(R^n\).
By \(w\)-coherence of \(R\), it is of \(w\)-finitely presented
type. Since $K/JK$ is of $w$-finite type and $L$ is of
$w$-finitely presented type over $R$,
Proposition~\ref{prop:finite-type-tools}(2), applied to the $w$-exact sequence
\[
0\longrightarrow T\longrightarrow K/JK\longrightarrow L\longrightarrow0,
\]
shows that $T$ is of $w$-finite type over $R$.
Thus $\Tor_1^S(R,X)$ is of $w$-finite type as well.
\end{proof}

\begin{corollary}
\label{cor:tor-free-lifting}
Let $M$ be a $w$-module over $R$, set $S=R\ltimes M$ and
$J=0\ltimes M$, and let $X$ be an $S$-module. Suppose that the
kernel of the multiplication map
\[
\mu_X:M\otimes_R(X/JX)\longrightarrow JX,\qquad
m\otimes\bar x\longmapsto(0,m)x,
\]
is GV-torsion over $R$. Then $X$ is of
$w$-finitely presented type over $S$ if and only if $X/JX$ is
of $w$-finitely presented type over $R$. In particular, this applies when $\mu_X$
is injective, as happens when
$X$ is flat over $S$ or when $X\cong S\otimes_R Y$ for an
$R$-module $Y$.
\end{corollary}

\begin{proof}
For a GV-torsion module $T$, the finitely generated submodule
$0\subseteq T$ satisfies $T/0=T$, so $T$ is of $w$-finite type. Theorem~\ref{thm:w-tor-lifting} therefore gives the
equivalence. By Lemma~\ref{lem:tor-multiplication}, injectivity
of $\mu_X$ is equivalent to vanishing of $\Tor_1^S(R,X)$.
For an induced module,
$S\otimes_R Y\cong Y\oplus(M\otimes_R Y)$ as an $R$-module,
and $\mu_X$ identifies $M\otimes_R Y$ with its second summand.
If $X$ is flat over $S$, the Tor group vanishes by flatness.
\end{proof}

\begin{remark}\label{rem:tor-hypothesis-scope}
For \(JX=0\), the sequence \(0\to J\to S\to R\to0\) gives
\(\Tor_1^S(R,X)\cong M\otimes_R X\).
Theorem~\ref{thm:inflated-fp} therefore proves both directions
of the Tor criterion for these modules without coherence of \(R\).
For arbitrary $X$, the converse can fail when $R$ is not
$w$-coherent, as Example~\ref{ex:tor-converse-fails} shows.
\end{remark}

\begin{proposition}
\label{prop:square-zero-patching}
Let $A$ be a ring, let $I$ be an ideal with $I^2=0$, put
$D=A/I$, and let $X$ be an $A$-module.
\begin{enumerate}[label=(\arabic*),leftmargin=*]
\item If $X/IX$ is finitely generated over $D$, then $X$ is
finitely generated over $A$.
\item If $X/IX$ is finitely presented over $D$ and
$\Tor_1^A(X,D)$ is finitely generated, then $X$ is finitely
presented over $A$.
\end{enumerate}
\end{proposition}

\begin{proof}
Part (1) is nilpotent Nakayama
\cite[Corollary~1.3.22]{wang-monograph2}. Indeed, lift generators
to a finitely generated submodule $Y\subseteq X$.
Then $W=X/Y$ satisfies $W=IW=I^2W=0$.
For (2), we use the standard Tor argument; cf.~\cite[Theorem~4.4.2]{glaz}.
By (1), choose an exact sequence
$0\to K\to A^n\to X\to0$ with $n$ finite. Tensoring with
$D$ gives
\[
\Tor_1^A(X,D)\longrightarrow K/IK\longrightarrow D^n
\longrightarrow X/IX\longrightarrow0.
\]
Put $L=\ker(D^n\to X/IX)$ and let
$\alpha:\Tor_1^A(X,D)\to K/IK$ be the map in this sequence.
Since $X/IX$ is finitely presented, $L$ is finitely generated
over $D$. By hypothesis, $\Tor_1^A(X,D)$ is finitely generated
over $D$, so its quotient
$\operatorname{im}\alpha\cong\Tor_1^A(X,D)/\ker\alpha$
is also finitely generated. Exactness gives
\[
0\longrightarrow\operatorname{im}\alpha\longrightarrow K/IK
\longrightarrow L\longrightarrow0.
\]
Thus $K/IK$ is finitely generated over $D$. Applying (1) to
$K$ shows that $K$ is finitely generated over $A$, and hence
$X$ is finitely presented over $A$.
\end{proof}

\begin{remark}
\label{rem:upper-superfluous}
Proposition~\ref{prop:square-zero-patching} applies to any ideal
$I\subseteq A$ with $I^2=0$, without a pullback construction.
The Tor hypothesis in part (2) holds, for example, if
$X$ is flat. It also holds if $D$ is Noetherian, $I$ is
finitely generated and $X$ is finitely generated. Tensoring
$0\to I\to A\to D\to0$ with $X$ gives the exact sequence
\[
0\longrightarrow\Tor_1^A(X,D)\longrightarrow X\otimes_A I
\xrightarrow{\mu}X\longrightarrow X/IX\longrightarrow0,
\qquad \mu(x\otimes i)=ix.
\]
Since $X$ and $I$ are finitely generated and $I^2=0$,
$X\otimes_A I$ is a finitely generated $D$-module.
As $D$ is Noetherian, its submodule $\Tor_1^A(X,D)$ is
also finitely generated.
\end{remark}

\Needspace{14\baselineskip}
\begin{remark}
\label{rem:tor-hypothesis-needed}
Finite generation of the square-zero ideal alone is insufficient, even in the
$w$-setting. Let $V$ be an infinite-dimensional $k$-vector space,
let $B=k\ltimes V$ with $V^2=0$, and put $I=0\ltimes V$.
Set $A=B\ltimes B$, $J=0\ltimes B$, and
$X=A/(0\ltimes I)$. The ideal $J$ is cyclic and $D=A/J=B$.
The rings $A$ and $B$ are zero-dimensional DW-rings, so their
$w$-finiteness notions are the ordinary ones. The cyclic module
$X$ satisfies $X/JX\cong B$, which is finitely presented over $D$.
However, $X$ is not finitely presented over $A$, since its relation
ideal $0\ltimes I$ is not finitely generated. Under the map
$\mu_X:B\otimes_B B\cong B\to JX\cong B/I$, multiplication
is the quotient map, so $\ker\mu_X\cong I$ is not finitely
generated. Thus $\Tor_1^A(D,X)\cong I$ is not of $w$-finite type.
This example shows that finite generation of $J$ cannot replace
the Tor hypothesis in Theorem~\ref{thm:w-tor-lifting}.
\end{remark}

For flat $M$,
the tensor-power formula is the specialization of
\cite[Section~4, observations~4.5(1) and Theorem~4.6(4)]{fgr}
to the exact functor $M\otimes_R-$.

\begin{proposition}
\label{prop:tor-periodicity}
Let $S=R\ltimes M$ and identify $R$ with $S/(0\ltimes M)$.
\begin{enumerate}[label=(\arabic*),leftmargin=*,nosep]
\item For every $M$, $\Tor_1^S(R,R)\cong M$ and
$\Tor_2^S(R,R)\cong M\otimes_R M$.
\item If $M$ is flat over $R$, then
\begin{equation}\label{eq:tor-periodicity}
\Tor_n^S(R,R)\cong M^{\otimes_R n}\qquad(n\ge1).
\end{equation}
\item If $M$ is a nonzero $w$-flat $w$-module of $w$-finite type, then
\[
\wfd_S R=\operatorname{fd}_S R=\wwgldim(S)=\infty.
\]
\end{enumerate}
\end{proposition}

\begin{proof}
\textup{(1)} Write $J=0\ltimes M$. The standard quotient calculation gives
$\Tor_1^S(R,R)=J/J^2\cong M$.
Dimension shifting gives
$\Tor_2^S(R,R)\cong\Tor_1^S(M,R)$.
Tensoring $0\to J\to S\to R\to0$ with the $S$-module $M$
identifies the latter group with $M\otimes_S J$, because
$JM=0$. This is $M\otimes_R M$.

\textup{(2)} For flat $M$, put $P_n=S\otimes_R M^{\otimes_R n}$,
where $M^{\otimes0}=R$, and define
\[
d_n(s\otimes m_1\otimes\cdots\otimes m_n)
=s(0,m_1)\otimes m_2\otimes\cdots\otimes m_n.
\]
As $R$-modules, $P_n=M^{\otimes n}\oplus M^{\otimes(n+1)}$,
and $d_n$ is $(a,b)\mapsto(0,a)$.
With the augmentation $S\to R$, this is an exact complex.
Each $P_n$ is flat over $S$ because it is induced from the
flat $R$-module $M^{\otimes n}$: for an $S$-module $Y$,
$Y\otimes_S P_n=Y\otimes_R M^{\otimes n}$.
It is therefore a flat resolution of $R$ over $S$.
Since $d_n(P_n)\subseteq JP_{n-1}$ and $JR=0$, tensoring with
$R$ gives
\[
R\otimes_S P_n\cong M^{\otimes_R n},\qquad
1_R\otimes d_n=0\quad(n\ge1).
\]
Consequently,
\[
\Tor_n^S(R,R)
\cong H_n(R\otimes_S P_\bullet)
\cong M^{\otimes_R n}\quad(n\ge1),
\]
which proves \eqref{eq:tor-periodicity}.

\textup{(3)} GV-torsion-freeness and $M\ne0$
give $\mathfrak m\in\wMax(R)$ with $M_{\mathfrak m}\ne0$.
This module is finitely generated and flat over the local
ring $R_{\mathfrak m}$, hence finite free of positive rank.
Applying the flat calculation at
$\mathfrak q=\mathfrak m\ltimes M\in\wMax(S)$ gives
\[
\Tor_n^{S_{\mathfrak q}}(R_{\mathfrak m},R_{\mathfrak m})
\cong (M_{\mathfrak m})^{\otimes_{R_{\mathfrak m}}n}\ne0
\quad(n\ge1).
\]
Thus $\operatorname{fd}_{S_{\mathfrak q}}R_{\mathfrak m}=\infty$.
The localization formulas for $w$-flat dimension and
$w$-weak global dimension, and
$\operatorname{fd}_S R\ge\wfd_S R$, give the result.
\end{proof}

\begin{corollary}
\label{cor:wwgd-infinite}
Let $M$ be a nonzero $w$-module of $w$-finite type that is
$w$-flat. Then
$\wwgldim(R\ltimes M)=\infty$.
\end{corollary}
\begin{proof}
This follows from Proposition~\ref{prop:tor-periodicity}(3).
\end{proof}

\begin{remark}
\label{rem:tor-flatness-needed}
The tensor-power formula is not a general periodicity
statement. If $R=k$ and $M=k^r$, its $n$th term has dimension
$r^n$, so these dimensions are not periodic for $r>1$.
Flatness is also needed in the calculation of higher Tor.
For $R=\mathbb Z$ and $M=\mathbb Z/2\mathbb Z$, write
$S=\mathbb Z[e]/(e^2,2e)$.
A free resolution of $R=S/(e)$ starts with maps
\[
S^5\xrightarrow{d_4}S^3\xrightarrow{d_3}S^2
\xrightarrow{d_2}S\xrightarrow{e}S\longrightarrow R\longrightarrow0,
\]
where
\[
d_2=(2\ \ e),\qquad
d_3=\begin{pmatrix}e&0&0\\0&2&e\end{pmatrix},\qquad
d_4=\begin{pmatrix}2&e&0&0&0\\0&0&e&0&0\\0&0&0&2&e\end{pmatrix}.
\]
These images equal the respective kernels, as follows by
writing each element of $S$ as $a+be$ with $a\in\mathbb Z$
and $b\in\mathbb F_2$. After tensoring with $R$, the third
homology is $(\mathbb Z/2\mathbb Z)^2$, whereas
$M^{\otimes_R3}\cong\mathbb Z/2\mathbb Z$.
\end{remark}

\begin{remark}
\label{rem:upper-corner}
Flatness of $T$ over $S$ in Theorem~\ref{thm:square-zero-square}
does not imply that $S$ has finite weak global dimension.
For $R=F=k$ and $M=k$, one has $T=S$ and
$\operatorname{fd}_S T=0$, but $S=k[e]/(e^2)$ has infinite weak
global dimension.
\end{remark}

\section{\texorpdfstring{$w$}{w}-properties of CPI-extensions}
\label{sec:cpi}

In this section $D$ is an integral domain with quotient field
$K$ and $P$ is a prime ideal. The classical CPI-extension of
Boisen and Sheldon \cite{boisen-sheldon} is
\[
C(D,P)=D+PD_P=\pi^{-1}(D/P),
\]
where $\pi:D_P\to D_P/PD_P=\operatorname{Frac}(D/P)$.
It fits into the square
\[
\begin{tikzcd}
C(D,P) \ar[r,hook] \ar[d,two heads] & D_P \ar[d,two heads]\\
D/P \ar[r,hook] & \operatorname{Frac}(D/P).
\end{tikzcd}
\]
The ordinary prime spectrum of $C(D,P)$ was determined by Boisen
and Sheldon \cite{boisen-sheldon}. We use this Milnor square to
describe the $w$-maximal ideals and compute the localizations at their
inverse images. We also apply general $w$-theoretic criteria for Milnor
squares to determine the $w$-coherence and strong Mori properties of
$C(D,P)$.
The results below apply to arbitrary prime ideals $P$.

When $P=0$ the ring is $D$. When $P$ is maximal the ring
is $D_P$: for $s\in D\setminus P$, write $as+p=1$ with
$a\in D$ and $p\in P$, giving $d/s=da+dp/s\in D+PD_P$.
The remaining case is $0\ne P$ with $P$ nonmaximal.

\begin{lemma}
\label{lem:conductor-w-ideal}
Suppose $0\ne P$ is nonmaximal, put $C=C(D,P)$,
$T=D_P$, and $\mathfrak c=PD_P$.
Then $(C:T)=\mathfrak c$ and $\mathfrak c_v=\mathfrak c$.
In particular $\mathfrak c$ is a proper $w$-ideal of $C$.
\end{lemma}

\begin{proof}
Since $0\ne P$ is nonmaximal, $D/P$ is a proper subring of its
quotient field and $PD_P$ is nonzero and maximal in $D_P$.
By \cite[Theorem~8.3.3(1)]{wang-monograph2}, we obtain
\[
(C:T)=\mathfrak c=\mathfrak c_v.
\]
Since $w\le v$, the conductor is also a $w$-ideal.
\end{proof}

\Needspace{24\baselineskip}
For the proposition below, assume that $0\ne P$ is a nonmaximal prime ideal
of $D$, and put $C=C(D,P)$, $\mathfrak c=PD_P$, and $B=D/P$.
Since $D\cap\mathfrak c=P$, the canonical quotient map induces an
isomorphism $C/\mathfrak c\cong B$. We write
$\rho:C\twoheadrightarrow B$, given explicitly by
$\rho(d+x)=d+P$ for $d\in D$ and $x\in\mathfrak c$; in particular,
$\ker\rho=\mathfrak c$.

\begin{proposition}
\label{prop:w-cpi-spectrum}
With this notation, the following assertions hold.
\begin{enumerate}[label=(\arabic*),leftmargin=*,nosep]
\item The $w$-maximal ideals of $C$ are given by
\begin{equation}\label{eq:wmax-cpi}
\wMax(C)=\{\rho^{-1}(\mathfrak n):\mathfrak n\in\wMax(B)\}.
\end{equation}
\item If $\mathfrak m/P=\mathfrak n\in\wMax(B)$ and
$\widetilde{\mathfrak m}=\rho^{-1}(\mathfrak n)$, then
\begin{equation}\label{eq:cpi-localization}
C_{\widetilde{\mathfrak m}}=D_{\mathfrak m}+PD_P
\quad\text{inside }D_P.
\end{equation}
\item The primes of $C$ below $\mathfrak c$ are the ideals
$\mathfrak pD_P$ for primes $\mathfrak p\subseteq P$ of $D$;
the primes above $\mathfrak c$ are the inverse images of
the primes of $B$. These two families meet at $\mathfrak c$.
\end{enumerate}
\end{proposition}

\begin{proof}
The ordinary spectrum description is the CPI spectrum
description of \cite{boisen-sheldon}.
One can see comparability directly: if $a\in C\setminus\mathfrak c$,
then $a$ is a unit of $D_P$, and $a^{-1}\mathfrak c=\mathfrak c$.
Thus $\mathfrak c\subseteq aC$, so every ideal of $C$ is
comparable with $\mathfrak c$.

Since $\mathfrak c=PD_P\ne0$ and
$B=D/P\subsetneq\operatorname{Frac}(B)$, the square is a Milnor
square of type I.
Since $B$ is not a field, choose a nonzero nonunit $b\in B$.
Then $bB$ is a nonzero proper $w$-ideal of $B$, and
\cite[Theorem~8.3.12]{wang-monograph2} implies that
$\rho^{-1}(bB)$ is a proper $w$-ideal of $C$ strictly containing
$\mathfrak c$.
As $\mathfrak c$ is $w$-closed by
Lemma~\ref{lem:conductor-w-ideal}, no $w$-maximal ideal lies
at or below $\mathfrak c$. The maximal $w$-ideal correspondence
\cite[Theorem~8.3.20]{wang-monograph2} therefore gives
\eqref{eq:wmax-cpi}.

For the localization formula, localizing $C=D+\mathfrak c$
at $D\setminus\mathfrak m$ gives
$D_{\mathfrak m}+\mathfrak c$, because every such denominator
is a unit of $D_P$ and preserves $\mathfrak c$.
If $u=d+j\in C\setminus\widetilde{\mathfrak m}$ with
$d\in D$ and $j\in\mathfrak c$, then $d\notin\mathfrak m$.
The factorization $u=d(1+j/d)$ shows that $u$ is invertible
in $D_{\mathfrak m}+\mathfrak c$: the inverse of
$1+j/d$ belongs to $1+\mathfrak c\subseteq C$.
This proves \eqref{eq:cpi-localization}.
\end{proof}

\begin{example}
\label{ex:cpi-localization}
Let $D=k[x,y]$, $P=(x)$, $C=C(D,P)$ and $\mathfrak m=(x,y)$.
Then $D/P=k[y]$ and $\mathfrak m/P=(y)$ is $w$-maximal.
Write $\widetilde{\mathfrak m}$ for its inverse image in $C$.
The element $x/y$ belongs to $PD_P\subseteq C$, but does
not belong to $D_{\mathfrak m}$. Thus
$C_{\widetilde{\mathfrak m}}\ne D_{\mathfrak m}$.
Moreover, in $C$ there is a strictly increasing chain
\[
xC\subsetneq (x/y)C\subsetneq (x/y^2)C\subsetneq\cdots.
\]
Strictness follows because equality at any step would imply
$1/y\in C$, whose residue $1/y$ is not in $k[y]$.
All these principal ideals are $w$-closed. Thus $C$ need not be
strong Mori even when $D$ is Noetherian and $D_P$ is a DVR.
\end{example}

\begin{theorem}
\label{thm:cpi-finiteness}
Let $P$ be a prime ideal of the domain $D$, and let
$C=C(D,P)$.
\begin{enumerate}[label=(\arabic*),leftmargin=*]
\item If $P=0$, then $C=D$.
\item If $P$ is maximal, then $C=D_P$.
\item If $0\ne P$ is nonmaximal, then $C$ is $w$-coherent
if and only if $D/P$ is $w$-coherent and $D_P$ is a valuation
domain. In this case $C$ is never a strong Mori domain.
\end{enumerate}
\end{theorem}

\begin{proof}
Parts (1) and (2) follow from $C(D,0)=D$ and
$C(D,P)=D_P$ for maximal $P$.
For (3), the common kernel $PD_P$ is nonzero and maximal in
the domain $D_P$, and $D/P$ is a proper subring of its
quotient field. Thus, by \cite[Theorem~8.5.22]{wang-monograph2},
$C$ is $w$-coherent if and only if $D/P$ is $w$-coherent and $D_P$ is a
valuation domain, since every valuation domain is $w$-coherent.
Suppose that $C$ were a strong Mori domain. By the strong Mori criterion
\cite[Theorem~8.4.16]{wang-monograph2}, $D/P$ would then have to be a
field. Since $P$ is nonmaximal, $D/P$ is not a field. Hence $C$ is not a
strong Mori domain.
\end{proof}

\begin{example}
\label{ex:cpi-degenerate}
The nonmaximality hypothesis in Theorem~\ref{thm:cpi-finiteness}(3) is essential. When $P$ is maximal, $C(D,P)=D_P$ can be strong Mori even though $D_P$ is not a valuation domain. For example, take $D=k[t^2,t^3]_{(t^2,t^3)}$ and
$P=(t^2,t^3)D$, the domain is one-dimensional Noetherian,
hence DW, so $P$ is a $w$-prime. Here $C=D_P=D$ is coherent
and strong Mori, but is not a valuation domain.
Indeed $t=t^3/t^2$ is integral over $D$ and not in $D$,
whereas valuation domains are integrally closed.
The condition $P\ne0$ is also necessary: for $D=\mathbb Z$ and
$P=(0)$, we have $C(D,P)=\mathbb Z$, which is a strong Mori domain.
Thus both $P\ne0$ and the nonmaximality of $P$ are necessary for the
conclusion in Theorem~\ref{thm:cpi-finiteness}(3) that $C(D,P)$ is not a
strong Mori domain.
\end{example}

To see why the classical Noetherian criterion
\cite[Theorem~8.4.15]{wang-monograph2} requires a finite residue field
extension, let $A=k+XK[[X]]$, where $k\subsetneq K$ is a field
extension. This ring is local with maximal ideal
$\mathfrak c=XK[[X]]$, and
$\mathfrak c/\mathfrak c^2\cong K$ as a $k$-vector space.
If $[K:k]=\infty$, then $A$ is not Noetherian, although
$k$ and $K[[X]]$ are Noetherian and $\mathfrak c$ is principal
over $K[[X]]$. If $[K:k]<\infty$, a $k$-basis of $K$
shows that $A$ is finitely generated over $k[[X]]$, so it is
Noetherian. Thus $A$ is Noetherian if and only if $[K:k]<\infty$.

\section*{Funding}
This work was partially supported by the National Natural Science Foundation
of China (Grant No.~11401493) and the Fundamental Research Funds for the
Central Universities, Southwest Minzu University (Grant No.~ZYN2023019).

\section*{Conflict of interest}
The author states that there is no conflict of interest.

\section*{Data availability}
No data was used for the research described in this article.

\end{document}